\documentclass[11pt,a4paper]{amsart}
\usepackage{amssymb}
\usepackage{amsmath}
\usepackage{amsfonts}
\usepackage{bbm}
\usepackage{enumerate}
\usepackage{amsthm}
\usepackage{amsbsy}
\numberwithin{equation}{section}

\newtheorem{theorem}{Theorem}[section]
\newtheorem{lemma}[theorem]{Lemma}

\newtheorem{corollary}[theorem]{Corollary}

\theoremstyle{definition}\newtheorem*{proof0}{Proof}
\newenvironment{pf}{\begin{proof0}
}{\mbox{}\hfill\qed\end{proof0}}

\theoremstyle{definition}\newtheorem{example}{Example}[section]

\theoremstyle{definition}\newtheorem{remark}{Remark}[section]

\begin{document}

\title{Polytope numbers and their properties}

\begin{abstract}
Polytope numbers for a polytope are a sequence of nonnegative
integers that are defined by the facial information of a polytope.
Every polygon is triangulable and a higher dimensional analogue of
this fact states that every polytope is triangulable, namely, every
polytope can be decomposed into simplexes. Thus it may be possible
to represent polytope numbers by sums of simplex numbers. We
analyzes a special type of triangulation, called pointed
triangulation, and develops several methods to represent polytope
numbers by sums of simplex numbers.
\end{abstract}

\maketitle \tableofcontents

\section*{Introduction}
\begin{figure}[h]
\begin{picture}(80,80)(0,0)
\put(-40,0){0} \put(-10,0){1} \put(-10,20){$\bullet$} \put(30,0){4}
\multiput(25,20)(10,0){2}{$\bullet$} \multiput(25,30)(10,0){2}{$\bullet$}
\multiput(27,22)(10,0){2}{\line(0,10){10}} \multiput(27,23)(0,10){2}{\line(10,0){10}}
\put(75,0){9} \multiput(65,20)(10,0){3}{$\bullet$}
\multiput(65,30)(10,0){3}{$\bullet$} \multiput(65,40)(10,0){3}{$\bullet$}
\multiput(67,22)(20,0){2}{\line(0,10){20}}
\multiput(67,23)(0,20){2}{\line(10,0){20}} \multiput(67,22)(10,0){2}{\line(0,10){10}}
\multiput(67,23)(0,10){2}{\line(10,0){10}} \put(128,0){16}
\multiput(115,20)(10,0){4}{$\bullet$} \multiput(115,30)(10,0){4}{$\bullet$}
\multiput(115,40)(10,0){4}{$\bullet$}
\multiput(115,50)(10,0){4}{$\bullet$} \multiput(117,22)(10,0){2}{\line(0,10){10}}
\multiput(117,23)(0,10){2}{\line(10,0){10}}
\multiput(117,22)(20,0){2}{\line(0,10){20}} \multiput(117,23)(0,20){2}{\line(10,0){20}}
\multiput(117,22)(30,0){2}{\line(0,10){30}}
\multiput(117,23)(0,30){2}{\line(10,0){30}}
\end{picture}
\begin{center}
\end{center}
\caption{Square numbers\label{figure0}}
\end{figure}
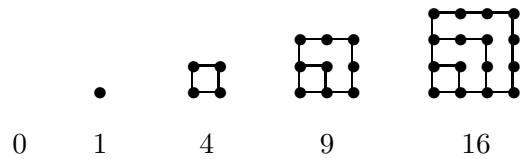

\emph{Polygonal numbers} are a sequence of nonnegative integers constructed
geometrically by a polygon. The square numbers are the numbers of
points in square arrays as in Figure \ref{figure0}.
\\ \emph{Polytope numbers} are higher dimensional analogues of polygonal numbers
or, equivalently, polygonal numbers are two dimensional polytope
numbers. Every polygon is triangulable and a higher dimensional
analogue of this fact is that every polytope is triangulable. Thus
it may be possible to represent polytope numbers by sums of simplex
numbers. We analyze a special type of triangulation, called pointed
triangulation, and develop several methods to represent polytope
numbers by sums of simplex numbers, which we formulate as decomposition theorems. We also consider several applications of polytope numbers to other mathematcial topics.

We provide basic definitions and notations in Section
\ref{section1}. We define pointed triangulation and consider
shellings for pointed triangulations in Section \ref{section2}. We
define polytope numbers and provide two descriptions of polytope
numbers in Section \ref{section3}. We furnish several ways of
decomposing polytope numbers into simplex numbers and formularize
them as decomposition theorems in Section \ref{section4}. We
illustrate decomposition theorems by applying them to several
polytopes in Section \ref{section5}. We suggest applications of
polytope numbers in Section \ref{section6}.

\section{Preliminaries\label{section1}}

Polytope numbers for a polytope, which we study in Section
\ref{section3}, are a sequence of natural numbers defined by a
recurrence relation that uses the facial information of a polytope.
Therefore properties of polytope numbers and those of polytopes
are closely related. We collect basic material for polytopes that we
use in our discussion. We begin with the definition of polytope and
provide several examples of polytopes. We next introduce the faces
of a polytope, the interior of a polytope, lines in general
position, linear functions in general position, polytopal complexes,
and shellings of polytopal complexes. The basis of this section is
the contents in Ziegler's book \cite{ziegler}.

\subsection{The definition of polytope}
Let $\mathbb{R}^d$ be the vector space of all column vectors of
length $d$ with real entries and $(\mathbb{R}^d)^*$ be its dual
vector space. The column vectors in $\mathbb{R}^d$ represent points.
The column vectors $\mathbf{0}$ and $\mathbf{1}$ are the column
vectors of all zeros and all ones, respectively, and the column
vectors $\mathbf{e}_1,\mathbf{e}_2,\ldots,\mathbf{e}_d$ are the unit
vectors in $\mathbb{R}^d$. Each row vector $\mathbf{a}\in
(\mathbb{R}^d)^*$ represents the linear form
$l_{\mathbf{a}}:\mathbb{R}^d\rightarrow \mathbb{R}$ defined by
$l_{\mathbf{a}}\mathbf{x}\mapsto \mathbf{ax}$. The row vectors
$\mathbbm{O}$ and $\mathbbm{1}$ denote the all zeros and all ones
row vectors in $(\mathbb{R}^d)^*$, respectively, and the row vectors
$\mathbf{f}_1,\mathbf{f}_2,\ldots,\mathbf{f}_d$ are the unit vectors
in $(\mathbb{R}^d)^*$. The vector $\mathbf{x}^{*}$ is the transpose
of $\mathbf{x}$.

The nonempty \emph{affine subspaces} are the translates of linear
subspaces. The set of all \emph{affine combinations} of a finite set
$\{\mathbf{x}_1,\mathbf{x}_2,\ldots,\mathbf{x}_n\}$ is
$$aff\big(\{\mathbf{x}_1,\mathbf{x}_2,\ldots,\mathbf{x}_n\}\big)\\
=\bigg\{\mathbf{x}\in\mathbb{R}^d\,\bigg|\,\mathbf{x}=\sum_{i\in[n]}
\lambda_i \mathbf{x}_i\ \text{for}\
\lambda_i\in\mathbb{R},\sum_{i\in[n]} \lambda_i=1\bigg\},$$
where $[n]=\{1,2,\ldots,n\}$. A set of
$n$ points is \emph{affinely independent} if its affine hull has
dimension $n-1$.

For a subset $K$ of $\mathbb{R}^d$, let $conv(K)$ be the convex hull
of $K$. A \emph{$\mathcal{V}$-polytope} is the convex hull of a
finite set of points in some $\mathbb{R}^d$. An
\emph{$\mathcal{H}$-polytope} is a bounded intersection of finitely
many closed halfspaces in some $\mathbb{R}^d$. Denoting by $A=(\mathbf{a}_1,\mathbf{a}_2,
\ldots,\mathbf{a}_m)^*$ an $m\times d$-matrix with the rows $\mathbf{a}_1,\mathbf{a}_2,
\ldots,\mathbf{a}_m$
and writing $\mathbf{z}=(z_1,z_2,\ldots,z_m)^*$, we can represent an 
$\mathcal{H}$-polytope by
$$P(A,\,\mathbf{z})=\big\{\mathbf{x}\in \mathbb{R}^d\,
\big|\,(A,\mathbf{z})\in\mathbb{R}^{m\times
d}\times \mathbb{R}^m,A\mathbf{x}\leq \mathbf{z}\big\}$$ where the inequality 
$A\mathbf{x}\leq
\mathbf{z}$ is the shorthand for a system of inequalities
$$\mathbf{a}_1 \mathbf{x}\leq z_1,\mathbf{a}_2\mathbf{x}\leq
z_2,\ldots,\mathbf{a}_m \mathbf{x}\leq z_m.$$ A
\emph{polytope} is a point set that is either a
$\mathcal{V}$-polytope or an $\mathcal{H}$-polytope.
\begin{theorem}[Main theorem for polytope \cite{ziegler}]
A subset of $\mathbb{R}^d$ is a $\mathcal{V}$-polytope if and only
if it is an $\mathcal{H}$-polytope.
\end{theorem}
\noindent The dimension of a polytope $P$, denoted by $dim(P)$, is
the dimension of its affine hull and a \emph{$d$-polytope} is a
polytope of dimension $d$ in some $\mathbb{R}^e$ with $e\ge d$.

Some recycling operations produce new polytopes. Let $P$ be a
$d$-polytope and $\mathbf{x_0}$ be a point outside of $aff(P)$ (for
this we embed $P$ into $\mathbb{R}^n$ for some $n>d$). A
\emph{pyramid} over $P$ is
$$pyr(P)=conv\big(P\cup \{\mathbf{x_0}\}\big).$$
The face set of $pyr(P)$ is
$$\big\{F,pyr(F)\,\big|\,\text{$F$ is a face of $P$}\big\}.$$
Similarly, a \emph{bipyramid} over $P$ is
$$bipyr(P)=conv\big(P\cup\{\mathbf{x}_{-},\mathbf{x}_{+}\}\big)$$
where both $\mathbf{x}_{-}$ and $\mathbf{x}_{+}$ are in outside of
$aff(P)$ and an interior point of the segment
$conv\big(\{\mathbf{x}_{-},\mathbf{x}_{+}\}\big)$ is an interior
point of $P$.

For two polytopes $P$ and $P^{\prime}$ the \emph{product} of $P$ and $P^{\prime}$ is
$$P\times P^{\prime}=\big\{(\mathbf{x},\mathbf{x}^{\prime})^*\, \big|\,
\mathbf{x}\in P,\mathbf{x}^{\prime}\in P^{\prime}\big\}.$$ The dimension of $P\times
P^{\prime}$ is $dim(P)+dim(P^{\prime})$ and the face set of $P\times P^{\prime}$ is
$$\big\{F\times F^{\prime}\,\big|\,\text{$F$ (resp. $F^{\prime}$) is a face of $P$ (resp. $P^{\prime}$).}\big\}$$

\begin{example}
A \emph{$d$-simplex} $\alpha^d$ is the convex hull of $d+1$ affinely
independent points in $\mathbb{R}^n$ with $n\ge d$. Thus a $d$-simplex
is a polytope of dimension $d$ with $d+1$ vertexes and it is a
pyramid over a $(d-1)$-simplex. The \emph{standard $d$-simplex}
$\alpha^d_s$ is the simplex in $\mathbb{R}^{d+1}$ defined by
$$\alpha_s^d=conv\big(\{\mathbf{e}_1,\mathbf{e}_2,\ldots,\mathbf{e}_{d+1}\}\big).$$

We construct a $d$-cross polytope $\beta^d$ by an iteration. Let
$\beta^0$ be a point. For $d\ge 1$ we define a \emph{$d$-cross
polytope} to be $\beta^d=bipyr(\beta^{d-1})$. The \emph{standard
$d$-cross polytope} $\beta^d_s$ is the cross polytope
defined by
$$\beta_s^d=conv\big(\{\pm\mathbf{e}_1,\pm\mathbf{e}_2,\ldots,\pm\mathbf{e}_d\}\big).$$

We also form a \emph{$d$-measure polytope} $\gamma^d$ by an
iteration. Let $\gamma^0$ be a point and $\gamma^1$ be a line
segment. For $d\ge 2$ we define
$\gamma^d=\gamma^{d-1}\times\gamma^1$. In particular, the
\emph{standard $d$-measure polytope} $\gamma^d_s$ is the
measure polytope defined by
$$\gamma_s^d=conv\bigg(\bigg\{\sum_{i\in[d]}a_i\mathbf{e}_i\,\bigg|\, a_i\in\{1,-1\}\bigg\}\bigg).$$
\end{example}

\subsection{The faces and the interior of a polytope}
For a polytope $P$ we define $vert(P)$ to be the vertex set of $P$
and $\mathcal{F}(P)$ (resp. $\mathcal{F}_k(P)$) to be the
face (resp. $k$-face) set of $P$.

For $\mathbf{y}\in P$ if every proper face of $P$ does not contain
$\mathbf{y}$, then we say that $\mathbf{y}$ is an \emph{interior
point} of $P$. We define $int(P)$ to be the \emph{interior} of $P$,
which is the set of all interior points in $P$, and
$\partial{P}=P\setminus int(P)$ to be the \emph{boundary} of $P$. We call
$relint(P)$ the \emph{relative interior} of P, which is the interior
of $P$ with respect to an embedding of $P$ into its affine hull where $P$ is full dimensional. Analogous to the interior of a
polytope, $relint(P)$ is the set of points in $P$ that are in no
proper face of $P$. By the definition of relative interior,
$$P=\underset{F\in \mathcal{F}(P)}{\biguplus
} relint(F)$$ where $\biguplus$ denotes the disjoint union.

\subsection{Lines and linear functions in general position} Let
$P=P(A,\mathbf{1})$ be a $d$-polytope with $A=(\mathbf{a}_1,\mathbf{a}_2,\ldots,\mathbf{a}_n)^*$. A line through
$\mathbf{0}\in int(P)$ is in \emph{general position} with respect to
$P$ if it is not parallel to any hyperplane that defines faces of
$P$ and it does not hit the intersection of any two of them. If we
write the line in the form $L(\mathbf{u})=\{t\,\mathbf{u}\, |\,
t\in\mathbb{R}\}$ for some $\mathbf{u}\neq\mathbf{0}$, then general
position means that $\mathbf{a}_i \mathbf{u}\neq \mathbf{a}_j
\mathbf{u}$ when $i\neq j$ and $\{i,j\}\subseteq[n]$.
\begin{lemma}[Ziegler \cite{ziegler}]\label{lgenpo0}
Let $P=P(A,\mathbf{1})$ and
$\mathbf{u}\in\mathbb{R}^d\setminus\{\mathbf{0}\}$. If $\lambda$ is
small enough, then the line $L(\mathbf{u}^{\lambda})$ is in general
position with respect to $P$ where
$$\mathbf{u}^{\lambda}=\mathbf{u}+(\lambda,\lambda^2,\ldots,\lambda^d)^*.$$
\end{lemma}
\begin{corollary}[Ziegler \cite{ziegler}]
For a polytope $P$ a line in general position with respect to $P$
exists.
\end{corollary}

A linear function $\mathbf{cx}$ is in general position with respect
to a polytope $P$ if it separates the vertexes of $P$, that is, if
$\mathbf{c}\mathbf{v}_i\neq\mathbf{c}\mathbf{v}_j$ for any two
distinct vertexes $\mathbf{v}_i$ and $\mathbf{v}_j$ of $P$.

\begin{lemma}[Ziegler \cite{ziegler}]\label{gpos0}
Let $P=P(A,\mathbf{1})$ and
$\mathbf{c}\in(\mathbf{R}^d)^*\setminus\{\mathbbm{O}\}$. If
$\lambda>0$ is small enough, then the linear function
$\mathbf{c}^{\lambda}\mathbf{x}$ is in general position with respect
to $P$ where
$$\mathbf{c}^{\lambda}=\mathbf{c}+(\lambda,\lambda^2,\ldots,\lambda^d).$$
\end{lemma}
\begin{corollary}[Ziegler \cite{ziegler}]\label{gpos1}
For a polytope $P$ a linear function in general position with
respect to $P$ exists.
\end{corollary}

Let $P$ be a polytope with $\mathbf{0}\in int(P)$ and $F$ be a
proper face of $P$. For a point $\mathbf{y}\in\mathbbm{R}^d$ we say
that $\mathbf{y}$ is a point \emph{beyond} $F$ if $\mathbf{y}$ and
$\mathbf{0}$ lie on different sides of $H_1$ for every facet
defining hyperplane $H_1$ that includes $F$, but on the same side of
$H_2$ for every facet defining hyperplane $H_2$ that does not
include $F$.

Let $P$ be a polytope. The vertexes and edges of $P$ form a finite,
undirected, and simple graph $G(P)$, called the \emph{graph} of $P$.
For every face $F$ of $P$ we define $G(F)$ to be the induced
subgraph of $G(P)$ on $vert(F)$, that is, the graph of all vertexes
in $F$ and all edges of $P$ between them. The graph $G(F)$ coincides
with the graph of $F$ if we consider $F$ as a polytope.

We consider an orientation of $G(P)$ that assigns a direction to
every edge. An orientation is \emph{acyclic} if no directed cycle is
in it. Thus, if a graph $G(P)$ has an acyclic orientation, then
$G(P)$ has a \emph{sink}, that is, a vertex that does not have an
edge directed away from it. If a linear function
$\mathbf{c}\mathbf{x}$ is in general position, then this linear
function gives a well-defined method to direct the graph $G(P)$, by
directing an edge $conv\big(\{\mathbf{v}_i,\mathbf{v}_j\}\big)$ from
$\mathbf{v}_i$ to $\mathbf{v}_j$ if
$\mathbf{c}\mathbf{v}_i>\mathbf{c}\mathbf{v}_j$. We call this
orientation the \emph{orientation} of $G(P)$ induced by
$\mathbf{c}\mathbf{x}$. Monotone paths on $P$, that is, edge paths
such that the objective function increases strictly in each step,
translate into directed paths in the orientation of $G(P)$ induced
by $\mathbf{c}\mathbf{x}$.

\begin{theorem}[Ziegler \cite{ziegler}]\label{usink}
For a polytope $P$, if $\mathbf{cx}$ is a linear function in general
position for $P$, then the orientation of $G(P)$ induced by
$\mathbf{c}\mathbf{x}$ is acyclic with a unique sink. This sink is
the unique point in $P$ such that $\mathbf{cx}$ achieves its minimum.
\end{theorem}
\noindent Note that if $O$ is an acyclic orientation of $G(P)$, then
the restriction of $G(P)$ to each nonempty subset $A$ of $vert(P)$
has a sink with respect to $O$.

An acyclic orientation $O$ of $G(P)$ is \emph{good} if for every
nonempty face $F$ of $P$ the graph $G(F)$ has exactly one sink. The
existence of good acyclic orientations of $G(P)$ follows from
Theorem \ref{usink}.
\begin{corollary}[Ziegler \cite{ziegler}]\label{gpos2} For a polytope
$P$, if $\mathbf{cx}$ is in general position for $P$, then it is in
general position for each face of $P$.
\end{corollary}

\subsection{Polytopal complexes and shellings} A \emph{polytopal
complex} $\mathcal{C}$ is a finite collection of polytopes that
satisfies the following conditions:
\begin{enumerate}[\emph{Condition} 1.]
\item \label{ccondition1}The empty polytope is in $\mathcal{C}$.

\item \label{ccondition2}If $P$ is an
element of $\mathcal{C}$, then every face of $P$ is also in
$\mathcal{C}$.

\item \label{ccondition3} The intersection
of two polytopes $P$ and $P^{\prime}$ in $\mathcal{C}$ is a face both of the
polytopes $P$ and $P^{\prime}$.
\end{enumerate}
The dimension of $\mathcal{C}$ is defined by
$$dim(\mathcal{C})=max\big\{dim(P)\,\big|\,P\in \mathcal{C}\big\}$$
and the \emph{underlying set} of $\mathcal{C}$ is the point set
$|\mathcal{C}|=\bigcup_{P\in\mathcal{C}}P$.

A polytopal complex is \emph{pure} if each of its faces is in a face
of dimension $dim(\mathcal{C})$, that is, if all of the inclusion
maximal faces of $\mathcal{C}$, called the \emph{facets} of
$\mathcal{C}$, have the same dimension. Let $\mathcal{C}_P$ be the
complex formed by the faces of a polytope $P$. The \emph{boundary
complex} $\mathcal{C}_{\partial P}$ is the subcomplex of
$\mathcal{C}_P$ formed by all proper faces of $\mathcal{C}_P$. Thus
its underlying set is
$$|\mathcal{C}_{\partial P}|=\partial P=P\backslash relint(P).$$ A
\emph{subdivision} of a polytope $P$ is a polytopal complex
$\mathcal{C}_P$ with the underlying space $|\mathcal{C}_P|=P$. The
subdivision is a \emph{triangulation} if every polytope in
$\mathcal{C}_P$ is a simplex.

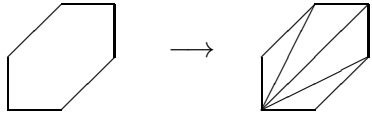
\begin{figure}[h]
\begin{center}
\begin{picture}(270,60)(-70,0)
\put(0,0){\line(1,0){20}} \put(0,0){\line(0,1){20}}
\put(20,0){\line(1,1){20}}\put(0,20){\line(1,1){20}}
\put(20,40){\line(1,0){20}}\put(40,20){\line(0,1){20}}
\put(60,20){$\longrightarrow$} \put(95,0){\line(1,0){20}}
\put(95,0){\line(0,1){20}}
\put(115,0){\line(1,1){20}}\put(95,20){\line(1,1){20}}
\put(115,40){\line(1,0){20}}\put(135,20){\line(0,1){20}}
\put(95,0){\line(1,2){20}}\put(95,0){\line(2,1){40}}
\put(95,0){\line(1,1){40}}
\end{picture}
\end{center}
\caption{A hexagon's triangulation\label{figure1}}
\end{figure}

Although a polytopal complex is a set of polytopes, we use it as a
generalization of a polytope and we define the faces of a polytopal
complex to be its elements. For example, the faces of a hexagon's
triangulation in Figure \ref{figure1} are the vertexes, the edges,
and the triangles in this triangulation.
\\
From this generalization, we denote
$$\mathcal{C}=\bigoplus_{F\in\mathcal{C}}F$$ and define
$\mathcal{F}_k(\mathcal{C})$ to be the set of $k$-faces of
$\mathcal{C}$. Since every face of a polytopal complex $\mathcal{C}$
is in a facet of $\mathcal{C}$, if $dim(\mathcal{C})=d$ and
$\mathcal{C}$ is a pure polytopal complex, then we denote
$$\mathcal{C}=\bigoplus_{F\in\mathcal{F}_d(\mathcal{C})}F.$$
Note that
$$\mathcal{C}=\bigoplus_{F\in\mathcal{C}}F\ \text{and}\
\mathcal{C}=\bigoplus_{F\in\mathcal{F}_d(\mathcal{C})}F$$ are
polytopal complex analogues of
$$|\mathcal{C}|=\bigcup_{F\in\mathcal{C}}F\ \text{and}\
|\mathcal{C}|=\bigcup_{F\in\mathcal{F}_d(\mathcal{C})}F,$$
respectively.

Let $\mathcal{C}_P$ be a polytopal complex obeying
$|\mathcal{C}_P|=P$. For $F\in\mathcal{F}(P)$ the subdivision
$\mathcal{C}_{F}$ of $F$ on $\mathcal{C}_{P}$ is
$$\mathcal{C}_{F}=\big\{F^{\prime}\in\mathcal{C}_{P}\,\big|\,F^{\prime}\subseteq F.\big\}.$$
Therefore
$$\mathcal{C}_{F}=\underset{F^{\prime}\in\mathcal{C}_{P}\atop
F^{\prime}\subseteq F}{\bigoplus} F^{\prime}.$$ We call
$\mathcal{C}_{F}$ the subdivision of $F$ with respect to
$\mathcal{C}_{P}$.

For a pure $d$-polytopal complex $\mathcal{C}$ a \emph{shelling} of
$\mathcal{C}$ is a linear ordering of the facets
$F_1,F_2,\ldots,F_s$ in $\mathcal{C}$ such that either $\mathcal{C}$
is $0$-dimensional or it satisfies the following conditions:
\begin{enumerate}[\emph{Condition} 1.]
\item \label{shellingcondition1}The boundary complex $\mathcal{C}_{\partial F_1}$
of the first facet $F_1$ has a shelling. 

\item\label{shellingcondition2}For $j\in[s]\setminus\{1\}$ the intersection of the
facet $F_j$ with the union of previous facets is a nonempty set and
a beginning segment of a shelling of the $(d-1)$-dimensional
boundary complex in $F_j$, that is,
$$F_j\cap\bigg(\bigcup_{i\in[j-1]}F_i\bigg)=F^{\prime}_1\cup F^{\prime}_2\cup\cdots\cup F^{\prime}_r$$
for some shelling $F^{\prime}_1,F^{\prime}_2,\ldots,F^{\prime}_r,\ldots,F^{\prime}_{s^{\prime}}$ of
$\mathcal{C}_{\partial F_j}$. In particular, this condition requires
that $F_j\cap\Big(\underset{i\in[j-1]}{\bigcup}F_i\Big)$ has a shelling,
therefore it has to be pure $(d-1)$-dimensional and connected for
$j\in[s]\setminus\{1\}$.
\end{enumerate}
A polytopal complex is \emph{shellable} if it is pure and has a
shelling.

\begin{remark}\mbox{}\hfill
\begin{enumerate}[1.]
\item Every simplex is shellable and every ordering of its facets
is a shelling. These facts immediately follow by induction on
dimension since the intersection of $F_j$ with $F_i$ for $i<j$ is
always a facet of $F_j$. \item For triangulations, the shelling
condition \ref{shellingcondition1} is redundant. Thus we can
simplify the shelling condition \ref{shellingcondition2}
considerably. We can replace the shelling condition
\ref{shellingcondition2} with
\begin{enumerate}[\emph{Condition} $2^{\prime}$.]
\item For $j\in[s]\setminus\{1\}$ the intersection of the facet $F_j$ with the
previous facets is nonempty and pure $(d-1)$-dimensional. In other
words, for every $i<j$ there exists some $l<j$ such that $F_i\cap
F_j$ is a subpolytope of $F_l\cap F_j$ and $F_l\cap F_j$ is a facet
of $F_j$.
\end{enumerate}
\end{enumerate}
\end{remark}

Let $P$ be a polytope and $\mathbf{x}$ be a point. The point
$\mathbf{x}$ lies in \emph{general position} with respect to the
polytope $P$ if $\mathbf{x}$ is not in the affine hull of a facet in
$P$. A facet $F$ of a $P$ is \emph{visible} from $\mathbf{x}$ if for
every $\mathbf{y}\in F$ the line segment
$conv\big(\{\mathbf{x},\mathbf{y}\}\big)$ intersects $P$ only in the
point $\mathbf{y}$. Equivalently, $F$ is visible from $\mathbf{x}$
if and only if $\mathbf{x}$ and $int(P)$ are on different sides of
the hyperplane $aff(F)$ spanned by $F$. For example, if
$\mathbf{x}_G$ is beyond the face $G$, then the facets that include
$G$ are exactly those that are visible from $\mathbf{x}_G$.

\begin{theorem}[Bruggersser \cite{bruggersser}]\label{shelling2} Let $P$ be a
polytope and $\mathbf{x}$ be a point outside $P$. If point
$\mathbf{x}$ lies in \emph{general position} with respect to $P$,
then the boundary complex $\mathcal{C}_{\partial P}$ of $P$ has a
shelling such that the facets of $P$ that are visible from
$\mathbf{x}$ come first.
\end{theorem}
\noindent This theorem shows that every polytope is shellable. In
Section \ref{section2}, we need to use Ziegler's proof of Theorem
\ref{shelling2} \cite{ziegler}, thus we provide that proof here.
\begin{pf} Let $\mathbf{x}$ be a point that lies in general position
with respect to $P$. We choose a line $l$ through both $\mathbf{x}$
and a point in general position for $P$. The properties we need are
that $l$ contains $\mathbf{x}$, $l$ hits the interior of $P$, and
any two different facet-defining hyperplanes $H$ and $H^{\prime}$ of
$P$ satisfy $l\cap H\neq l\cap H^{\prime}$. For simplicity, we
assume that $l$ is not parallel to any of the facet hyperplanes,
thus there is no intersection point at infinity. We orient $l$ from
$P$ to x.

Now we imagine that $P$ is a polyhedral planet and there is a rocket
that starts on its surface at the point where the oriented line $l$
leaves the planet. This point lies on a unique facet $F_1$ of $P$
and for the first few minutes of the flight only $F_1$ is visible
from the rocket.

After a while, a new facet appears on the horizon. The rocket passes
through a hyperplane $H_2$ and we label the corresponding facet
$F_2$. We continue to label the facets $F_3,F_4,\ldots$ of $P$ in
the order where the rocket passes through their hyperplanes, that
is, in the order such that the facets appear on the horizon, becoming
visible from the rocket. Now we imagine that the rocket passes
through and comes back to the planet from the opposite side. We
continue the shelling by taking the facets in the order such that the
rocket passes though the hyperplanes $aff(F_i)$, that is, the
corresponding facets disappear on the horizon.

This rocket flight clearly gives a well-defined ordering on the
facets of $P$. What's more, the facets that are visible from
$\mathbf{x}$ form a beginning segment, since we see exactly those
facets at the point where the rocket passes through $\mathbf{x}$.

To verify that the ordering is a shelling, we consider the
intersection
$$\partial F_j\cap\bigg(\bigcup_{i\in[j-1]}F_i\bigg).$$
If we add $F_j$ before we pass through infinity, then this
intersection is exactly the set of those facets in $F_j$ that are
visible from the point $l\cap aff(F_j)$ where $F_j$ appears on the
horizon. Thus we know by induction on dimension that this collection
of facets of $F_j$ is shellable and can be continued to a shelling
of the whole boundary $\partial F_j$.

After the rocket passes through infinity, the intersection is the
family of nonvisible facets. This family of nonvisible facets is
also shellable because reversing the orientation of $l$ yields the
shelling with the reversed ordering of the facets.
\end{pf}

\begin{corollary}[Ziegler \cite{ziegler}]\label{arshelling}
For any two facets $F$ and $F^{\prime}$ of a polytope $P$, there is
a shelling of $\partial P$ such that $F$ is the first facet and
$F^{\prime}$ is the last one.
\end{corollary}
\begin{pf}
We choose a shelling line $l$ that intersects the boundary of $P$ in
the facets $F$ and $F^{\prime}$. For example, we choose two points
$\mathbf{x}$ and $\mathbf{x}^{\prime}$ beyond $F$ and $F^{\prime}$,
respectively, and let $l$ be the line determined by $\mathbf{x}$ and
$\mathbf{x}^{\prime}$. If necessary, we perturb $l$ to general
position.
\end{pf}

\section{Pointed triangulations \label{section2}}
Every polygon is triangulable and we may regard a polytope's
triangulation as a higher dimensional analogue of a polygon's
triangulatoin. We introduce a special kind of triangulation, called
pointed triangulation, and study shellings of this triangulation.
Pointed triangulations and their shellings constitute the main tool
to relate polytopes and polytope numbers.

Let $\mathcal{C}_{P}$ be a triangulation of a $d$-polytope $P$ that
satisfies the following conditions: For $d\in\mathbb{N}=\{0,1,\ldots\}$ let $[d]_0=\{0,1,\ldots,d\}$.
\begin{enumerate}[\emph{Condition} 1.]
\item \label{tcondition1}For $k\in[d]_0$ each $k$-face $F$ of $P$ has a triangulation
$$\mathcal{C}_{F}=\underset{\alpha^k\in\mathcal{C}_{F}}{\bigoplus}\alpha^k$$
such that there is a designated vertex $\mathbf{v}_{F}\in vert(F)$,
called the \emph{apex} of $F$, satisfying $$\mathbf{v}_{F}\in
\underset{\alpha^k\in\mathcal{C}_{F}}{\bigcap}\alpha^k.$$
\item\label{tcondition2} For any two faces $F_1$ and $F_2$ of $P$ if
$\{\mathbf{v}_{F_1},\mathbf{v}_{F_2}\}\subseteq F_1\cap F_2$ then
$\mathbf{v}_{F_1}=\mathbf{v}_{F_2}$.
\item\label{tcondition3} For each face $F$ of $P$ if $\mathbf{w}\in vert(F)\setminus\{\mathbf{v}_{F}\}$ then the
edge $conv\big(\{\mathbf{v}_{F},\mathbf{w}\}\big)$ is
in $\mathcal{C}_{F}$.
\end{enumerate}

We define
$$V(P)=\big\{\mathbf{v}_{F}\,\big|\,F\in\mathcal{F}(P)\big\}$$ and call
it the \emph{set of apexes} of $\mathcal{C}_P$. Note that $V(P)$ is
a multiset by the condition \ref{tcondition2}. Even though the set
$V(P)$ is dependent on a pointed triangulation of $P$, we use $V(P)$
by abusing notation. We call $\mathcal{C}_P$ the \emph{$V(P)$-pointed
triangulation}. Conditions
\ref{tcondition1}--\ref{tcondition3} are called the \emph{pointed
triangulation conditions}.
\begin{theorem}\label{dec1} Every polytope has a pointed triangulation.
\end{theorem}
\begin{pf} Let $P$ be a $d$-polytope. By
Lemma \ref{gpos1} a linear function in general position
$\mathbf{cx}$ with respect to $P$ exists and by Corollary
\ref{gpos2} the linear function $\mathbf{cx}$ decides a unique sink
$\mathbf{v}_{F}$ for each face $F\in\mathcal{F}(P)$. For each
$F\in\mathcal{F}(P)$ we define
$$V(F)=\big\{\mathbf{v}_G\,|\,G\in\mathcal{F}(F)\big\}$$
and
\begin{align*}
\mathcal{C}_F=&\{\emptyset\}\cup\Big\{conv\big(\{\mathbf{v}_{G_1},
\mathbf{v}_{G_2},\ldots,\mathbf{v}_{G_{k}}\}\big)\,
 \Big|\,k\in[d],\\
 &\hspace{3em}G_i\in\mathcal{F}(F),G_i\supsetneq G_{i+1},
 \mathbf{v}_{G_i}\notin G_{i+1}\Big\}.
\end{align*}
We claim that $\mathcal{C}_{P}$ is the $V(P)$-pointed triangulation.
To verify this claim, we use induction on
dimension.

We first show that $\mathcal{C}_P$ is a triangulation of $P$. By the
definition of $\mathcal{C}_P$ the polytopal complex $\mathcal{C}_P$
is a triangulation, thus we need only show that $|\mathcal{C}_P|=P$.
Moreover, $|\mathcal{C}_P|\subseteq P$, therefore it suffices to show that
$|\mathcal{C}_P|\supseteq P$.

Let $\mathbf{x}\in P$. Since
$$P=\bigcup_{F\in\mathcal{F}_{d-1}(P)\atop \mathbf{v}_P\notin F}
conv\big(\{\mathbf{v}_P\}\cup F\big),$$ there is a facet $F$ of $P$
such that
$$\begin{cases}
\mathbf{v}_P\notin F\\
\mathbf{x}\in conv\big(\{\mathbf{v}_P\}\cup F\big)\end{cases}.$$ By
the way, for each $F\in\mathcal{C}$ the complex $\mathcal{C}_F$
satisfies $|\mathcal{C}_F|=F$ by the induction hypothesis, hence there is
a $(d-1)$-simplex $\alpha^{d-1}$ in $\mathcal{C}_F$ such that
$$\mathbf{x}\in conv\big(\{\mathbf{v}_P\}\cup
\alpha^{d-1}\big).$$ By the definition of $\mathcal{C}_P$ the
polytope $conv\big(\{\mathbf{v}_P\}\cup \alpha^{d-1}\big)$ is an
element of $\mathcal{C}_P$, thus $\mathbf{x}\in|\mathcal{C}_P|$. It
follows that $|\mathcal{C}_P|\supseteq P$.

We now show that $\mathcal{C}_P$ is the $V(P)$-pointed triangulation.
By the definition of $\mathcal{C}_P$, the
complex $\mathcal{C}_P$ satisfies the pointed triangulation
conditions \ref{tcondition1} and \ref{tcondition3}. Let $F_1$ and
$F_2$ be two faces of $P$ such that
$\{\mathbf{v}_{F_1},\mathbf{v}_{F_2}\}\subseteq F_1\cap F_2$. Since
$$\mathbf{v}_{F_1}\in F_1\cap F_2\in\mathcal{F}(F_1),$$ we obtain
$\mathbf{v}_{F_1\cap F_2}=\mathbf{v}_{F_1}$, and similarly,
$\mathbf{v}_{F_1\cap F_2}=\mathbf{v}_{F_2}$. This yields
$\mathbf{v}_{F_1}=\mathbf{v}_{F_2}$. As a result, $\mathcal{C}_{P}$
satisfies the pointed triangulation condition \ref{tcondition2}.
This proves that $C_{P}$ is the $V(P)$-pointed triangulation.
\end{pf}

From now on, every pointed triangulation is formed by the method in
the proof of Theorem \ref{dec1}.

\begin{theorem}\label{decshell}
Let $P$ be a $d$-polytope with the $V(P)$-pointed triangulation
$\mathcal{C}_{P}$. There is a shelling of
$\mathcal{C}_{P}$ such that for
$\mathcal{F}_d(\mathcal{C}_{P})=\big\{\alpha^d_1,\alpha^d_2,\ldots,\alpha^d_s\big\}$
if $j\in[s]\setminus\{1\}$, then the number of facets of $\alpha^d_j$ in
$\alpha^{d_j}\cap\bigg(\underset{i\in[j-1]}{\bigcup}\alpha^d_i\bigg)$, denoted
by $l_j$, satisfies $l_j\in[q-1]$.
\end{theorem}
\begin{pf} The complex $\mathcal{C}_{P}$ is a triangulation,
thus we need only show the shelling condition $2^{\prime}$. We
borrow the definitions and notations in the proof of Theorem
\ref{shelling2}.

Let
$$\mathcal{F}_{d-1}(\mathcal{C}_{P},\widehat{\mathbf{v}_{P}})
=\big\{\alpha^{d-1}_1,\alpha^{d-1}_2,\ldots,\alpha^{d-1}_s\big\}$$
be the set of $(d-1)$-simplexes $\alpha^{d-1}$ on $\partial(P)$
satisfying
$$\begin{cases}
\mathbf{v}_{P}\notin\alpha^{d-1}\\
\alpha^d_i=conv\big(\{\mathbf{v}_{P}\}\cup\alpha^{d-1}_i\big)\end{cases}.$$
Since $\mathcal{C}_P$ is a pointed triangulation, we may instead
prove that there is a shelling of
$\mathcal{F}_{d-1}(\mathcal{C}_{P},\widehat{\mathbf{v}_{P}})$ such
that the number of facets of $\alpha^{d-1}_j$ in\\
$\alpha^{d-1}_j\cap\bigg(\underset{i\in[j-1]}{\bigcup}\alpha^{d-1}_i\bigg)$,
denoted by $l_j$, satisfies $l_j\in[d-1]$. We define such a
shelling to be the \emph{triangulation shelling}.

Without loss of generality, we assume that $\mathbf{f}_d\mathbf{x}$
is a linear function in general position that determines $V(P)$. We
choose a point $\mathbf{x}$ that lies in general position, that is,
does not lie in the affine hull of a facet of $P$, and a line $l$
through $\mathbf{x}$ and a point in general position. By Corollary
\ref{arshelling}, we may assume that $l$ passes through two distinct
facets $F^{d-1}$ and $F^{d-1}_{P}$ of $P$ such that $\mathbf{v}_m\in
F^{d-1}$ and $\mathbf{v}_{P}\in F^{d-1}_{P}$ where
$$\begin{cases}
\mathbf{f}_d\mathbf{v}_m=max\big\{\mathbf{f}_d\mathbf{v}\,\big|\,\mathbf{v}\in
vert(P)\big\}\\
\mathbf{f}_d\mathbf{v}_{P}=min\big\{\mathbf{f}_d\mathbf{v}\,\big|\
\mathbf{v}\in vert(P)\big\}\end{cases}.$$ We may further assume by
Lemma \ref{lgenpo0} that the line $l$ is orthogonal to the plane
$\mathbf{f}_d\mathbf{x}=0$ and the rocket defined by $l$ moves in the
direction that the value of $d$th coordinate in its position vector
increases. As in Theorem \ref{shelling2}, this rocket flight assigns
a shelling to the facets of $P$, thus it assigns a shelling to those
facets that do not contain $\mathbf{v}_{P}$.

If either $d=0$ or $d=1$, then the line $l$ evidently assigns a triangulation
shelling to
$\mathcal{F}_{d-1}(\mathcal{C}_{P},\widehat{\mathbf{v}_{P}})$.
Suppose that $d\ge 2$. We assume that the line $l$ assigns a
triangulation shelling to every pointed triangulation of a
$k$-polytope when $k\in[d-1]\cup\{0\}$.

We first assign a shelling to the facets of $P$ that do not contain
$\mathbf{v}_{P}$ as in Theorem \ref{shelling2}. Let them be
$F^{d-1}_1,F^{d-1}_2,\ldots$ and $\mathcal{C}_{F^{d-1}_i}$ be the
$V(F^{d-1}_i)$-pointed triangulation. In this case, the
indexes in the facets of $P$ equal those in the shelling of $P$. By
the induction hypothesis, $\mathcal{C}_{F^{d-1}_1}$ has a
triangulation shelling. Suppose that
$\underset{i\in[k]}{\bigcup}\mathcal{C}_{F^{d-1}_i}$ has a triangulation
shelling and let $\mathbf{v}_{F^{d-1}_{k+1}}$ be the sink of
$F^{d-1}_{k+1}$. Then $$\mathbf{v}_{F^{d-1}_{k+1}}\notin
\underset{i\in[k]}{\bigcup}F^{d-1}_i\ \text{or}\
\mathbf{v}_{F^{d-1}_{k+1}}\in\underset{i\in[k]}{\bigcup}F^{d-1}_i.$$

We suppose that $\mathbf{v}_{F^{d-1}_{k+1}}\notin\underset{i\in[k]}{\bigcup}
F^{d-1}_i$ and define
$$C_{F^{d-1}_{k+1}}=\underset{i\in[m_{k+1}]}{\bigoplus}\alpha^{d-1}_{k+1,i}.$$
If we use induction on dimension by considering the formation of
$\underset{i\in[k]}{\bigcup}\mathcal{C}_{F^{d-1}_i}$, then the triangulation
$\mathcal{C}_{F^{d-1}_{k+1}}$ has a shelling such that the number of
facets of $\alpha^{d-1}_{k+1,i}$ in
$$\bigg|\bigcup_{i\in[k]}\mathcal{C}_{F^{d-1}_i}\bigg|
\bigcup\bigg(\bigcup_{j\in[i]}\alpha^{d-1}_{k+1,j}\bigg)$$ is at most
$d-1$.

Suppose that $\mathbf{v}_{F^{d-1}_{k+1}}\in\underset
{i\in[k]}{\bigcup}F^{d-1}_i$. If every facet
$F_{k+1,1}^{d-2}$ of $F^{d-1}_{k+1}$ obeying
$$dim\bigg(\bigg|\underset{i\in[k]}{\bigcup}\mathcal{C}_{F^{d-1}_i}\bigg|\bigcap
F_{k+1,1}^{d-2}\bigg)=d-2$$ contains the apex
$\mathbf{v}_{F^{d-1}_{k+1}}$, then
$\underset{i\in[k+1]}{\bigcup}\mathcal{C}_{F^{d-1}_i}$ has a triangulation
shelling. Assume that $F^{d-1}_{k+1}$ has another facet
$F_{k+1,2}^{d-2}$ satisfying
$$\begin{cases}
F_{k+1,2}^{d-2}\subseteq\underset{i\in[k]}{\bigcup}F^{d-1}_i\\
dim\Big(F_{k+1,1}^{d-2}\bigcap F_{k+1,2}^{d-2}\Big)=d-3\\
\mathbf{v}_{F^{d-1}_{k+1}}\notin F_{k+1,2}^{d-2}\end{cases}.$$ Using
Corollary \ref{arshelling}, we can choose a shelling of
$\partial(F^{d-1}_{k+1})$ such that $F_{k+1,1}^{d-2}$ is the last
facet and $F_{k+1,2}^{d-2}$ is the first one. If we use induction on
dimension by considering the formation of
$\underset{i\in[k]}{\bigcup}\mathcal{C}_{F^{d-1}_i}$, then we can assign a
triangulation shelling to
$\underset{i\in[k+1]}{\bigcup}\mathcal{C}_{F^{d-1}_i}$.

Considering all of the facets in
$\mathcal{F}_{d-1}(\mathcal{C}_{P},\widehat{\mathbf{v}_{P}})$ allows us to assign a triangulation shelling to these facets.
\end{pf}
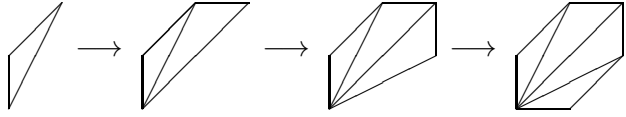
\begin{figure}[h]
\begin{center}
\begin{picture}(250,60)(-20,-5)
\put(0,0){\line(0,1){20}}\put(0,20){\line(1,1){20}}
\put(0,0){\line(1,2){20}} \put(25,20){$\longrightarrow$}

\put(50,0){\line(0,1){20}}
\put(50,20){\line(1,1){20}}\put(50,20){\line(1,1){20}}
\put(70,40){\line(1,0){20}} \put(50,0){\line(1,2){20}}
\put(50,0){\line(1,1){40}}\put(95,20){$\longrightarrow$}

\put(120,0){\line(0,1){20}} \put(120,20){\line(1,1){20}}
\put(120,0){\line(1,2){20}}\put(140,40){\line(1,0){20}}
\put(120,0){\line(1,1){40}} \put(120,0){\line(2,1){40}}
\put(160,20){\line(0,1){20}}\put(165,20){$\longrightarrow$}

\put(190,0){\line(0,1){20}} \put(190,20){\line(1,1){20}}
\put(190,0){\line(1,2){20}}\put(210,40){\line(1,0){20}}
\put(190,0){\line(1,1){40}} \put(190,0){\line(1,0){20}}
\put(190,0){\line(2,1){40}}\put(210,0){\line(1,1){20}}
\put(230,20){\line(0,1){20}}
\end{picture}
\end{center}
\caption{A triangulation shelling of a pointed triangulation}
\end{figure}

\section{Polytope numbers\label{section3}} We define
polytope numbers in this section. By using the definition of
polytope numbers, we also derive the product formula for polytope
numbers. Intuitively, polytope numbers for a polytoope are a
sequence of numbers associated with a polytope. To count these
numbers more effectively, we need a canonical method to describe
them geometrically, called the geometric description of polytope
numbers. We begin with the geometric description of simplex numbers,
introduce the facet-cut, suggest the geometric description of
polytope numbers, and finally consider the description of polytope
numbers by vertex sets.

\subsection{The definition of polytope numbers} For a $d$-polytope
we define a sequence of \emph{polytpe numbers} $P(n)$ and
\emph{interior polytope numbers} $P(n)^\sharp$ by double induction
on $d$ and $n$.

When $d=0$, we define $$\begin{cases}P(0)=0,P(n)=1\ \text{for}\ n\ge 1\\
P(0)^\sharp=0,P(n)^\sharp=1\ \text{for}\ n\ge 1\end{cases}.$$

When $d\ge 1$, we suppose that for each $k$-polytope $F$ satisfying
$k<d$ the numbers $F(n)$ and $F(n)^{\sharp}$ are defined. We define
$$\begin{cases}
P(0)=0,P(1)=1\\
P(n)=P(n-1)+\underset{\mathbf{v}_{P}\notin F\in\mathcal{F}(P)}{\sum}
F(n)^\sharp\ \text{for}\ n\ge 2\end{cases}$$ and $$\begin{cases}
P(0)^\sharp=0,P(1)^\sharp=0\\
P(n)^\sharp=P(n)-\underset{F\in\mathcal{F}(P)\backslash\{P\}}{\sum}
F(n)^\sharp\ \text{for}\ n\ge 2\end{cases},$$
where $\mathbf{v}_{P}$ is a fixed
vertex of $P$. Following the definition of polytope numbers, we can
also define
$$P(n)=\sum_{F\in\mathcal{F}(P)}F(n)^{\sharp}.$$
\begin{remark}\mbox{}\hfill
\begin{enumerate}[1.]
\item To define polytope numbers for a polytope $P$, we need to choose
a vertex $\mathbf{v}_{P}$ of $P$ and a vertex $\mathbf{v}_{F}$ of
$F$ for each $F\in\mathcal{F}(P)\setminus\{P\}$ with
$\mathbf{v}_{P}\notin F$. Defining $V(P)$ to be the set of such
vertexes, we call $P(n)$ the $V(P)$-polytope numbers.
\item The set $V(P)$ for the $V(P)$-polytope numbers coincides with the set
$V(P)$ for the $V(P)$-pointed triangulation defined in Section \ref{section2}.
The geometric description of polytope numbers explains the reason
for this coincidence.
\end{enumerate}
\end{remark}

\subsection{Polytope numbers for products of polytopes} Following
the definition of polytope numbers and that of the product of
polytopes, we can expect that polytope and interior polytope numbers
for the product of two polytopes are the product of two
corresponding polytope and interior polytope numbers, respectively. We shot that this is actually true.

For $i\in[2]$ let $P_i$ be a $d_i$-polytope. Suppose that
$P_i(n)$ are the $V(P_i)$-polytope numbers and $P_1\times P_2(n)$
are the $V(P_1\times P_2)$-polytope numbers. We claim that
\begin{equation}\label{pproduct}
\begin{cases}
P_1\times P_2(n)=P_1(n)\times P_2(n)\\
P_1\times P_2(n)^{\sharp}=P_1(n)^{\sharp}\times
P_2(n)^{\sharp}\end{cases}.\end{equation} To prove this claim, we
use double induction on the numbers $d=d_1+d_2$ and $n$.

If either $d=0$ or $n=0$, then the identity (\ref{pproduct}) is
trivially true. Suppose that $d\ge 1$ and $n\ge 1$. The apex of $P_1\times P_2$ is $\mathbf{v}_{P_1}\times \mathbf{v}_{P_2}$, thus 
\begin{align}\label{product}P_1\times P_2(n)&=P_1\times
P_2(n-1)+\sum_{\mathbf{v}_{P_1}\times \mathbf{v}_{P_2}\notin F_1\times F_2\atop
F_1\times F_2\in \mathcal{F}(P_1\times P_2)}F_1\times
F_2(n)^{\sharp}
\end{align}
By the induction hypothesis $$\begin{cases}P_1\times
P_2(n-1)=P_1(n-1)\times P_2(n-1)\\
F_1\times F_2(n)^{\sharp}=F_1(n)^{\sharp}\times
F_2(n)^{\sharp}\end{cases},$$
hence the identity (\ref{product}) becomes
\begin{align*}
P_1\times P_2(n)&=P_1(n-1)\times
P_2(n-1)+\sum_{\mathbf{v}_{P_1}\times\mathbf{v}_{P_2}\notin F_1\times F_2\atop
F_1\times F_2\in \mathcal{F}(P_1\times P_2)}F_1(n)^{\sharp}\times
F_2(n)^{\sharp}
\end{align*}
In addition,
\begin{align*}
P_1(n-1)\times P_2(n-1)&=
\underset{\mathbf{v}_{P_1}\in F_1\atop F_1\in\mathcal{F}(P_1)}{\sum}F_1(n)^{\sharp}
\times \underset{\mathbf{v}_{P_2}\in F_2\atop F_2\in\mathcal{F}(P_2)}{\sum}F_1(n)^{\sharp}\\
&=\sum_{\mathbf{v}_{P_1}\times \mathbf{v}_{P_2}\in F_1\times F_2\atop
F_1\times F_2\in \mathcal{F}(P_1\times P_2)}F_1(n)^{\sharp}\times
F_2(n)^{\sharp},
\end{align*}
therefore
\begin{align*}
P_1\times P_2(n)&=\sum_{
F_1\times F_2\in \mathcal{F}(P_1\times P_2)}F_1(n)^{\sharp}\times
F_2(n)^{\sharp}\\
&=\bigg(\sum_{F_1\in\mathcal{F}(P_1)}F_1(n)^{\sharp}\bigg)\times
\bigg(\sum_{F_2\in\mathcal{F}(P_2)}F_2(n)^{\sharp}\bigg)
=P_1(n)\times P_2(n).
\end{align*}
It follows that
$$P_1\times P_2(n)=P_1(n)\times P_2(n).$$
Similarly, $$P_1\times P_2(n)^{\sharp}=P_1(n)^{\sharp}\times
P_2(n)^{\sharp}.$$

Generalizing the identity (\ref{pproduct}) to the product of several
polytopes, we can compute polytope numbers for products of several
polytopes.
\begin{theorem}[Polytope numbers for the product of polytopes]
Let $P_1,P_2,\ldots,P_l$ be polytopes. Suppose that $P_i(n)$ are the $V(P_i)$-polytope numbers. If
$\underset{i\in[l]}{\prod}P_i(n)$ are the $V\bigg(\underset{i\in[l]}{\prod}P_i\bigg)$-polytope numbers, then
$$\begin{cases}
\bigg(\underset{i\in[l]}{\prod}P_i\bigg)(n)=\underset{i\in[l]}{\prod}P_i(n)\\
\bigg(\underset{i\in[l]}{\prod}P_i\bigg)(n)^{\sharp}=\underset{i\in[l]}{\prod}P_i(n)^{\sharp}
\end{cases}.$$
\end{theorem}

\subsection{The geometric description of simplex numbers and the facet-cut}
\label{geo} To
define polytope numbers for a polytope $P$, we need to use its
facial information. Hence it may be possible to arrange points in
$P$ that correspond to polytope numbers for $P$. Since polytope
numbers for $P$ are determined by a pointed triangulation
$\mathcal{C}_P$, we need to consider a method to describe simplex
numbers by such arrangements of points in simplexes. We proceed by
induction.

We first consider polytope numbers for the standard $d$-simplex
$\alpha^d_s$. We claim that
\begin{equation}\label{salpha}
\alpha_s^d(n)=\Big|\big\{\mathbf{x}\in\mathbb{N}^{d+1}\,\big|\,
\mathbbm{1}\mathbf{x}=n-1\big\}\Big|.
\end{equation}
If either $d=0$ or $(d,n)\in\big(\mathbb{N}\setminus\{0\}\big)\times [n]_1$,
then the identity (\ref{salpha}) is obviously true. For $d\ge 2$
and $n\ge 2$, we suppose that
$$\left\{\begin{array}{ll}
\alpha_s^k(n)=\Big|\big\{\mathbf{x}\in\mathbb{N}^{k+1}\,\big|\,
\mathbbm{1}\mathbf{x}=n-1\big\}\Big|&
\text{when}\ k\in[d-1]\\
\alpha_s^d(m)=\Big|\big\{\mathbf{x}\in\mathbbm{N}^{d+1}\,\big|\,
\mathbbm{1}\mathbf{x}= m-1\big\}\Big| &
\text{when}\ m\in[n-1]\end{array}\right..$$ Since
\begin{align*}
&\Big|\big\{\mathbf{x}\in\mathbb{N}^{d+1}\,\big|\,
\mathbbm{1}\mathbf{x}=n-1\big\}\Big|\\
=&\Big|\big\{\mathbf{x}\in\mathbb{N}^{d+1}\,\big|\,
\mathbbm{1}\mathbf{x}=n-2\big\}\Big|
=\Big|\big\{\mathbf{x}\in\mathbb{N}^{d}\,\big|\,
\mathbbm{1}\mathbf{x}=n-1\big\}\Big|\\
=&\alpha_s^d(n-1)+\alpha_s^{d-1}(n)=\alpha^d(n-1)+\alpha^{d-1}(n)\\
=&\alpha^d(n),
\end{align*}
this proves the claim.
\begin{lemma} For each $n$ the number $\alpha_s^d(n)$
is the number of points in
$$S_{\alpha_s^d}(n)=\big\{\mathbf{x}\in\mathbb{N}^{d+1}\,\big|\,
\mathbbm{1}\mathbf{x}=n-1\big\}.$$
\end{lemma}

We now consider the case of ordinary simplexes. Let $\alpha^d$ be a
$d$-simplex with
$vert(\alpha^d)=\{\mathbf{v}_1,\mathbf{v}_2,\ldots,\mathbf{v}_{d+1}\}$.
The sets $vert(\alpha_s^d)$ and $vert(\alpha^d)$ are affinely
independent, thus there is a bijective affine map $f_{\alpha^d}$
such that
$$f_{\alpha^d}(\mathbf{e}_i)=\mathbf{v}_i$$
when $i\in[d+1]$. It follows that $f_{\alpha^d}(\alpha_s^d)=\alpha^d$. Since
$f_{\alpha^d}$ is a bijection, we can correlate $\alpha^d(n)$ to the
points in $f_{\alpha^d}\big(S_{\alpha_s^d}(n)\big)$.

Let $\alpha^d$ be a $d$-simplex and $F$ be a facet of a $d$-simplex
$conv\Big(f_{\alpha^d}\big(S_{\alpha_s^d}(n)\big)\Big)$. Then the
number of points in $F\cap f_{\alpha^d}\big(S_{\alpha_s^d}(n)\big)$
is $$\alpha^{d-1}(n)={n+d-1-1\choose d-1}.$$ Now eliminating all
points of $F\cap f_{\alpha^d}\big(S_{\alpha_s^d}(n)\big)$ from
$f_{\alpha^d}\big(S_{\alpha_s^d}(n)\big)$ changes
$f_{\alpha^d}\big(S_{\alpha_s^d}(n)\big)$ into one of
$$S^i_{\alpha^d}(n)=f_{\alpha^d}\big(S^i_{\alpha_s^d}(n)\big)$$
for $i\in[d+1]$ where
$$\begin{cases}
S^i_{\alpha_s^d}(n)=\big\{\mathbf{x}\in\mathbb{N}^{d+1}\,\big|\,
\mathbbm{1}\mathbf{x}= n-1,\mathbf{f}_{i}\mathbf{x}\ge
1\big\}\\
\Big|S^i_{\alpha_s^d}(n)\Big|={n-1+d-1\choose d}
\end{cases}.$$
We call this process the \emph{facet-cut} and we represent it by
$$\alpha^d(n)-\alpha^{d-1}(n)=\alpha^d(n-1).$$ In general, successive $k$
facet-cuts on $\alpha^d(n)$ yield $\alpha^d(n-k)$.

\subsection{The geometric description of polytope numbers} We have geometrically
described simplex numbers by arranging points in a simplex. By using
this description, we consider the geometric description of polytope
numbers for a polytope by arranging points in a polytope.

Let $P$ be a $d$-polytope with the $V(P)$-pointed triangulation
$\mathcal{C}_{P}$. Assuming that $P(n)$ are
the $V(P)$-polytope numbers, we define
$$\begin{cases}
S_P(n)=\underset{\alpha^d\in\mathcal{F}_d(C_P)}{\bigcup}S_{\alpha^d}(n)\\
S_{P}(n)^\sharp=S_P(n)-\underset{\alpha^{d-1}\in\mathcal{F}_{d-1}(C_P)
\atop\alpha^{d-1}\subseteq\partial(P)}{\bigcup}S_{\alpha^{d-1}}(n)\end{cases}.$$
By using double induction on $d$ and $n$, we show that
$$\begin{cases}
P(n)=\big|S_{P}(n)\big|\\
P(n)^{\sharp}=\big|S_{P}(n)^{\sharp}\big|\end{cases}.$$

Suppose that $d=0$. By definition,
$$\begin{cases}
\big|S_{P}(0)\big|=0,\big|S_{P}(n)\big|=1\ \text{for}\ n\ge 1\\
\big|S_{P}(0)^\sharp\big|=0,\big|S_{P}(n)^\sharp\big|=1\ \text{for}\ n\ge 1
\end{cases}.$$ Assume that $d\ge 1$ and $n\in[1]_0$. Similarly,
$$\begin{cases}
\big|S_{P}(0)\big|=1,\big|S_{P}(1)\big|=1\\
\big|S_{P}(0)^\sharp|=0,|S_{P}(1)^\sharp\big|=0\end{cases}.$$ From
now on, we assume that $d\ge 2$ and $n\ge 2$.

For a point $p$ and a set $S$ we define $p+S=\{p+s\,|\, s\in S\}$.
By the definition of $S_P(n)$,
\begin{align*}
S_P(n)&=\big(\mathbf{v}_P+S_P(n-1)\big)\biguplus\Big\{S_P(n)\setminus\big(\mathbf{v}_P+S_P(n-1)\big)\Big\}\\
      &=\big(\mathbf{v}_P+S_P(n-1)\big)\biguplus
      \bigg(\bigcup_{{\mathbf{v}_P\notin\alpha^{d-1}\atop\alpha^{d-1}\in
      \mathcal{F}_{d-1}(\mathcal{C}_P)}}S_{\alpha^{d-1}}(n)\bigg).
      \end{align*}
Since $\mathcal{C}_P$ is a pointed triangulation, by induction on $d$
$$\bigcup_{{\mathbf{v}_P\notin\alpha^{d-1}\atop\alpha^{d-1}\in\mathcal{F}_{d-1}
      (\mathcal{C}_P)}}S_{\alpha^{d-1}}(n)=
      \bigcup_{{\mathbf{v}_P\notin F\atop F\in\mathcal{F}(P)}}
      S_F(n)=\biguplus_{{\mathbf{v}_P\notin F\atop F\in\mathcal{F}(P)}}S_F(n)^{\sharp},$$
which yields
$$S_P(n)=\big(\mathbf{v}_P+S_P(n-1)\big)\biguplus\bigg(\biguplus_{{\mathbf{v}_{P}\notin F\atop F\in\mathcal{F}(P)}}
S_F(n)^{\sharp}\bigg).$$
      Therefore, by double induction on $d$ and $n$,
\begin{align*}
\big|S_P(n)\big|&=\bigg|\big(\mathbf{v}_P+S_P(n-1)\big)\biguplus\bigg(\biguplus_{
{\mathbf{v}_P\notin
      F\atop F\in\mathcal{F}(P)}}S_F(n)^{\sharp}\bigg)\bigg|\\
      &=\big|S_P(n-1)\big|+\sum_{{\mathbf{v}_P\notin F\atop F\in\mathcal{F}(P)}}
      \big|S_F(n)^{\sharp}\big|=P(n-1)+\sum_{{\mathbf{v}_P\notin F\atop F\in\mathcal{F}(P)}}F(n)^{\sharp}\\
      &=P(n).
      \end{align*}
Similarly,
\begin{align*}
|S_{P}(n)^\sharp|
&=\bigg|S_P(n)\setminus\bigcup_{\alpha^{d-1}\in\mathcal{F}_{d-1}(C_P)
 \atop\alpha^{d-1}\subseteq\partial(P)}S_{\alpha^{d-1}}(n)\bigg|\\
&=\bigg|S_{P}(n)\setminus\bigcup_{F\in\mathcal{F}(P)\setminus\{P\}}S_{F}(n)\bigg|
 =\bigg|S_{P}(n)\setminus\biguplus_{F\in\mathcal{F}(P)\setminus\{P\}}S_{F}(n)^\sharp\bigg|\\
&=\big|S_{P}(n)\big|-\sum_{F\in\mathcal{F}(P)\setminus\{P\}}\big|S_F(n)^\sharp\big|
 =P(n)-\sum_{F\in\mathcal{F}(P)\setminus\{P\}}F(n)^\sharp\\
                              &=P(n)^\sharp.
                              \end{align*}
\begin{theorem}
Let $P$ be a $d$-polytope with the $V(P)$-pointed triangulation
$\mathcal{C}_P$. Suppose that $P(n)$ are the $V(P)$-polytope numbers.
Then \begin{align*}
P(n)&=\bigg|\bigcup_{\alpha^d\in\mathcal{F}_d(C_P)}S_{\alpha^d}(n)\bigg|,\\
P(n)^\sharp&=\Bigg|
\bigg(\bigcup_{\alpha^d\in\mathcal{F}_d(C_P)}S_{\alpha^d}(n)\bigg)
\setminus\bigg(\bigcup_{\alpha^{d-1}\in\mathcal{F}_{d-1}(C_P)
\atop\alpha^{d-1}\subseteq\partial(P)}S_{\alpha^{d-1}}(n)\bigg)\bigg|.
\end{align*}
\end{theorem}

\subsection{The vertex description of polytope numbers}
The $V(P)$-polytope numbers are determined by both the set $V(P)$ and the facial information of $P$. For example, for each $n$ the number
$\alpha^d_s(n)$ is, by the geometric description of polytope
numbers, the number of points in the set
$$S_{\alpha^d_0}(n)=\big\{\mathbf{x}\in\mathbb{N}^{d+1}\,\big|\,
\mathbbm{1}\mathbf{x}=n-1\big\}.$$ By the
way,
\begin{align*}\big\{\mathbf{x}\in\mathbb{N}^{d+1}\,\big|\,
\mathbbm{1}\mathbf{x}=n-1\big\}
&=\bigg\{\sum_{j\in[n-1]}\mathbf{e}_{i_j}\,\bigg|\,1\leq i_1\leq
i_2\leq\cdots\leq i_{n-1}\leq d+1\bigg\}\\
&=\bigg\{\sum_{i\in[n-1]}\mathbf{v}_{F_i}\,\bigg|\,
F_i\in\mathcal{F}\big(\alpha^d_s\big),F_i\supseteq
F_{i+1}\bigg\},\end{align*} hence for each $n$ the number $\alpha^d_s(n)$
equals the number of points in the set
$$\bigg\{\sum_{i\in[n-1]}\mathbf{v}_{F_i}\,\bigg|\,
F_i\in\mathcal{F}\big(\alpha^d_s\big),F_i\supseteq F_{i+1}\bigg\}.$$
In general, if $P$ is a $d$-polytope, then for each $n$ the number
$P(n)$ is the number of points in the set
$$\bigcup_{\alpha^d\in\mathcal{F}(\mathcal{C}_P)}f_{\alpha^d}\big(\alpha_s^d\big),$$
thus it may be possible to describe polytope numbers by vertex sets.

Let $P$ be a $d$-polytope with the $V(P)$-pointed triangulation
$\mathcal{C}_P$. For 
$n\ge 0$ we define two sequences of
sets $T_P(n)$ and $T_P(n)^{\sharp}$ by double induction
on the numbers $d$ and $n$.

Let $d=0$. We define
$$\begin{cases}T_P(0)=T_P(0)^{\sharp}=\emptyset\\
T_P(1)=T_P(1)^{\sharp}=\{\mathbf{0}\}\end{cases},$$ and for $n\ge
1$
$$T_P(n)=T_P(n)^{\sharp}=\{\mathbf{v}_P\}.$$

Let $d\ge 1$. We define
$$\begin{cases}
T_P(0)=T_P(0)^{\sharp}=\emptyset\\
T_P(1)=T_P(1)^{\sharp}=\{\mathbf{0}\}\end{cases}$$ and 
$$\begin{cases}
T_P(n)&=\underset{F\in\mathcal{F}(P)}{\bigcup}\big(\mathbf{v}_F+T_F(n-1)\big)
      =\bigg\{\underset{i\in[n-1]}{\sum}\mathbf{v}_{F_i}\,
\bigg|\,F_i\supseteq
F_{i+1}\bigg\},\\
T_P(n)^{\sharp}&=T_P(n)\setminus
\bigg(\underset{F\in\mathcal{F}(P)\setminus\{P\}}{\bigcup}T_F(n)\bigg)\\
&=\bigg\{\underset{i\in[n-1]}{\sum}\mathbf{v}_{F_i}\,\bigg|\,F_i\supseteq F_{i+1},
conv\big(\{\mathbf{v}_{F_1},\mathbf{v}_{F_2},\ldots,\mathbf{v}_{F_{n-1}}\}\big)
\subseteq
\partial(P)\bigg\}.
\end{cases}$$
for $n\ge
2$.

We claim that
\begin{equation}\label{vdescription}\begin{cases}
\big|T_P(n)\big|=P(n)\\
\big|T_P(n)^{\sharp}\big|=P(n)^{\sharp}\end{cases}.\end{equation} If
either $d=0$ or $n=0$, then the identities in (\ref{vdescription})
are apparent. Suppose that $d\ge 1$ and $n\ge 1$. Then
\begin{align*}
T_P(n)&=\bigcup_{F\in\mathcal{F}(P)}\big(\mathbf{v}_F+T_F(n-1)\big)\\
&=\big(\mathbf{v}_P+T_P(n-1)\big)\biguplus\bigg\{\bigcup_{
{\mathbf{v}_P\notin F\atop F\in\mathcal{F}(P)}}\big(\mathbf{v}_F+T_F(n-1)\big)\bigg\}\\
&=\big(\mathbf{v}_P+T_P(n-1)\big)\biguplus
\bigg\langle\bigcup_{{\mathbf{v}_P\notin F\atop
F\in\mathcal{F}(P)}}\bigg[T_F(n)^{\sharp}\bigcup\bigg\{
\bigcup_{{\mathbf{v}_F\in G\atop
G\in\mathcal{F}(F)\setminus\{F\}}}\big(\mathbf{v}_G+T_G(n-1)\big)\bigg\}\bigg]\bigg\rangle.
\end{align*}
Since
\begin{align*}
&\bigcup_{{\mathbf{v}_P\notin F\atop F\in\mathcal{F}(P)}}
\bigg[T_F(n)^{\sharp}\bigcup\bigg\{\bigcup_{{\mathbf{v}_F\in
G\atop G\in\mathcal{F}(F)\setminus\{F\}}}\big(\mathbf{v}_G+T_G(n-1)\big)\bigg\}\bigg]\\
=&\bigcup_{{\mathbf{v}_P\notin
F\atop F\in\mathcal{F}(P)}}T_F(n)^{\sharp}
=\biguplus_{{\mathbf{v}_P\notin F\atop
F\in\mathcal{F}(P)}}T_F(n)^{\sharp},
\end{align*}
we obtain $$T_P(n)=\big(\mathbf{v}_P+T_P(n-1)\big)\biguplus\bigg(
\biguplus_{{\mathbf{v}_P\notin F\atop
F\in\mathcal{F}(P)}}T_F(n)^{\sharp}\bigg).$$ Therefore, by double induction
on $d$ and $n$,
\begin{align*}
\big|T_P(n)\big|&=\bigg|\big(\mathbf{v}_P+T_P(n-1)\big)\biguplus\bigg(
\biguplus_{{\mathbf{v}_P\notin
F\atop F\in\mathcal{F}(P)}}T_F(n)^{\sharp}\bigg)\bigg|\\
&=\big|\mathbf{v}_P+T_P(n-1)\big|+\sum_{{\mathbf{v}_P\notin
F\atop F\in\mathcal{F}(P)}}\big|T_F(n)^{\sharp}\big|\\
&=P(n-1)+\sum_{{\mathbf{v}_P\notin F\atop F\in\mathcal{F}(P)}}F(n)^{\sharp}\\
&=P(n).
\end{align*}

By the definition of $T_P(n)^{\sharp}$ 
\begin{align*}
T_P(n)^{\sharp}&=T_P(n)\setminus
\bigg(\bigcup_{F\in\mathcal{F}(P)\setminus\{P\}}T_F(n)\bigg)\\
&=T_P(n)\setminus\bigg(\biguplus_{F\in\mathcal{F}(P)\setminus\{P\}}T_F(n)^{\sharp}\bigg),
\end{align*}
therefore
\begin{align*}
\big|T_P(n)^{\sharp}\big|&=\bigg|T_P(n)\setminus\bigg(\biguplus_{F\in\mathcal{F}(P)\setminus\{P\}}
T_F(n)^{\sharp}\bigg)\bigg|\\
&=\big|T_P(n)\big|-\sum_{F\in\mathcal{F}(P)\setminus\{P\}}\big|T_F(n)^{\sharp}\big|\\
&=P(n)-\sum_{F\in\mathcal{F}(P)\setminus\{P\}}F(n)^{\sharp}\\
&=P(n)^{\sharp}.
\end{align*}
\begin{theorem}
For a polytope $P$ if $P(n)$ are the $V(P)$-polytope numbers,
then
$$\begin{cases}
P(0)=P(0)^{\sharp}=\big|\emptyset\big|\\
P(1)=P(1)^{\sharp}=\big|\{\mathbf{0}\}\big|\end{cases}$$
and
$$\begin{cases}
P(n)=\bigg|\bigg\{\underset{i\in[n-1]}{\sum}\mathbf{v}_{F_i}\, \bigg|\,
F_i\supseteq F_{i+1}\bigg\}\bigg|,\\
P(n)^{\sharp}=\bigg|\bigg\{\underset{i\in[n-1]}{\sum}\mathbf{v}_{F_i}\,
\bigg|\, F_i\supseteq F_{i+1},conv\big(\{\mathbf{v}_{F_1},\mathbf{v}_{F_2},\ldots,\mathbf{v}_{F_{n-1}}\}\big)
\subseteq\partial(P)\bigg\}\bigg|.
\end{cases}$$
when $n\ge 2$.
\end{theorem}

\subsubsection{Computations of polytope numbers by the vertex description of polytope numbers}
Kim computed polytope numbers for
regular polytopes by the definition of polytope numbers
 \cite{kim}. We compute polytope numbers for several polytopes by the vertex
description of polytope numbers and show that our computation coincides Kim's computation.

For the standard $d$-simplex $\alpha^d_s$ we suppose that
$V\big(\alpha_s^d\big)$ is formed by the linear function
$$L_{\alpha^d}(x_1,x_2,\ldots,x_{d+1})=\sum_{i\in[d+1]}ix_i.$$
If $F$ is a face of $\alpha^d_s$, then
$$\begin{cases}
F=conv\{\mathbf{e}_{i_1},\mathbf{e}_{i_2},\ldots,\mathbf{e}_{i_k}\}\\
\mathbf{v}_F=\mathbf{e}_{i_1}\end{cases}$$ where $1\leq
i_1<i_2<\cdots<i_k\leq d+1$. Thus
$$T_{\alpha_s^d}(n)=\bigg\{(x_1,x_2,\ldots,x_{d+1})\in\mathbb{N}^{d+1}\,\bigg|\,
\sum_{i\in[d+1]}x_i=n-1 \bigg\},$$ which yields
$$\alpha^d(n)=\alpha_s^d(n)=\big|T_{\alpha^d}(n)\big|={n-1+d\choose d}.$$

For the standard $d$-cross polytope $\beta_s^d$ let
$$L_{\beta_s^d}(x_1,x_2,\ldots,x_d)=\sum_{i\in[d]}ix_i$$ be a linear function
that forms $V(\beta_s^d)$. If $F$ is a face of $\beta_s^d$,
then
$$F=conv\{a_{i_1}\mathbf{e}_{i_1},\ldots,a_{i_k}\mathbf{e}_{i_k}\}$$
where
$$\begin{cases}
\{i_1,i_2,\ldots,i_k\}\subseteq[d]\\
\{a_{i_1},a_{i_2},\ldots,a_{i_k}\}\subseteq\{1,-1\}\end{cases}.$$
Let $I_F^+$ and $I_F^-$ be subsets of $\{i_1,i_2,\ldots,i_k\}$
defined by
$$\begin{cases}
a_i=1&\text{for}\ i\in I_F^+\\
a_i=-1&\text{for}\ i\in I_F^-
\end{cases}.$$
For each face $F$ of $\beta_s^d$, the apex $\mathbf{v}_F$ is
$$\mathbf{v}_F=\begin{cases}
\mathbf{e}_1&\text{for}\ F=P\\
\mathbf{e}_{min(I_F^+)}&\text{for}\ I_F^+\neq\emptyset\\
-\mathbf{e}_{max(I_F^-)}&\text{for}\ I_F^+=\emptyset
\end{cases}.$$

We define a partial order $\prec$ on $V(\beta_s^d)$ by
$$\mathbf{v}_{F_1}\prec\mathbf{v}_{F_2}\ \text{for}\ F_1\subseteq F_2.$$
Then
$$\begin{cases}
\mathbf{e}_i\prec\mathbf{e}_j,-\mathbf{e}_i\prec-\mathbf{e}_j&\text{if}\
i\leq
j\\
\mathbf{e}_i\succ-\mathbf{e}_j&\text{if}\ i\neq j
\end{cases}.$$
Therefore $\beta^d(n)$ is the number of lattice points
$(x_1^{+}-x_1^{-},x_2,x_3,\ldots,x_d)$ satisfying
$$\begin{cases}
\ x_1^{+},x_1^{-}\ge 0\\
x_1^{+}+x_1^{-}+\underset{i\in[d]\setminus\{1\}}{\sum}|x_i|=n-1\end{cases}.$$
Recall that a lattice point is a point each of whose coordinates is an integer. We can
easily show that the number of such lattice points equals the number
of lattice points $(x_1,x_2,\ldots,x_d)$ obeying
$$\begin{cases}
x_1\ge 0\\
x_1+\underset{i\in[d]\setminus\{1\}}{\sum}|x_i|\leq n-1\end{cases}.$$

For the standard $d$-measure polytope $\gamma_s^d$ let
$$L_{\gamma_s^d}(x_1,x_2,\ldots,x_d)=-\bigg(\sum_{i\in[d]}(d+1-i)x_i\bigg)$$
be a linear function in general position for $\gamma_s^d$ that forms
$V\big(\gamma_s^d\big)$. A face $F$ of $\gamma^d$ is determined by
the intersection of the following hyperplanes that defines facets of
$\gamma_s^d$;
$$x_{i_1}=a_{i_1},x_{i_2}=a_{i_2},\ldots,x_{i_k}=a_{i_k}$$
where
$$\begin{cases}\{i_1,i_2,\ldots,i_k\}\subseteq[d]\\
\{a_{i_1},a_{i_2},\ldots,a_{i_k}\}\subseteq[1]_0\end{cases},$$ which gives 
$\mathbf{v}_F=\underset{j\in[k]}{\sum}a_{i_j}\mathbf{e}_{i_j}$. Thus
$$\gamma^d(n)=\gamma_s^d(n)=\Big|\big\{(x_1,x_2,\ldots,x_d)\big|\,x_i\in[n-1]_0\big\}\Big|.$$

Let $$\alpha_k^{d-1}=\bigg\{(x_1,x_2,\ldots,x_d)\, \bigg|\,
\sum_{i\in[d]}x_i=k,x_i\ge 0\bigg\}$$ be a hypersimplex and
$$L_{\alpha_k^{d-1}}(x_1,x_2,\ldots,x_d)=\sum_{i\in[d]}(d+1-i)x_i$$ be a
linear function in general position for $\alpha_k^{d-1}$ that forms
$V(\alpha_k^{d-1})$. Since the hyperplanes
$x_i=0$ and $x_i=1$ for $i\in[d]$ determine the facets of $\alpha_k^{d-1}$, the
hyperplanes
$$x_{i_1}=a_{i_1},x_{i_2}=a_{i_2},\ldots,x_{i_k}=a_{i_k}$$
satisfying
$$\begin{cases}\{i_1,i_2,\ldots,i_k\}\subseteq[d]\\
\{a_{i_1},a_{i_2},\ldots,a_{i_k}\}\subseteq[1]_0\end{cases},$$
determine a face $F$ of $\alpha^{d-1}_k$. Let
$$\begin{cases}I_F^0=\big\{i\in\{i_1,\ldots,i_k\}\,\big|\,
a_i=0\big\}\\
I_F^1=\big\{i\in\{i_1,i_2,\ldots,i_k\}\,\big|\,
a_i=1\big\}\end{cases}$$ and $I_F$ be the
set of the first $k-\big|I_F^1\big|$ minimal indexes in
$$[d]\setminus\big(I_F^0\cup I_F^1\big).$$ Note that $\big|I_F^0\big|\leq d-k$ and
$\big|I_F^1\big|\leq k$. Then
$$\mathbf{v}_F=\sum_{i\in
I_F^1}\mathbf{e}_i+\sum_{i\in I_F}\mathbf{e}_i.$$ 

We claim that $$\alpha_k^{d-1}(n)=\bigg|\bigg\{
(x_1,x_2,\ldots,x_d)\,\bigg|\,
\sum_{i\in[d]}x_i=k(n-1),x_i\in[n-1]_0\bigg\}\bigg|.$$ To
establish this claim, we need only show that for each
$(x_1,x_2,\ldots,x_n)$ with
$$\begin{cases}
\underset{i\in[d]}{\sum}x_i=k(n-1)\\
x_i\in[n-1]_0\end{cases}$$ there is a set
$\{\mathbf{v}_{F_1},\mathbf{v}_{F_2},\ldots,\mathbf{v}_{F_{n-1}}\}$
in $V(P)$ such that $$\begin{cases} F_i\supseteq F_{i+1}\\
\underset{i\in[n-1]}{\sum}\mathbf{v}_{F_i}=(x_1,x_2,\ldots,x_d)\end{cases}.$$

Let
$\mathbf{x}_1=(x_{1,1},x_{1,2},\ldots,x_{1,d})$ with
$$\begin{cases}
\underset{i\in[d]}{\sum} x_{1,i}=k(n-1)\\
x_{1,i}\in[n-1]_0
\end{cases}.$$ We define $I_1$ to be the set of
indexes $i$ such that $x_{1,i}=n-1$, $J_1$ to be the set of
$k-|I_1|$ smallest numbers from $[d]-I_1$, and 
$$\begin{cases}
\mathbf{y}_1=\underset{i\in I_1\cup J_1}{\sum}\mathbf{e}_i\\
\mathbf{x}_2=\mathbf{x}_1-\mathbf{y}_1\end{cases}.$$ Inductively, we
define $I_l$ to be the set of indexes $i$ such that $x_{l,i}=n-l$,
$J_l$ to be the set of $k-|I_l|$ smallest numbers from
$[d]-I_l$, and
$$\begin{cases}
\mathbf{y}_l=\underset{i\in I_l\cup J_l}{\sum}\mathbf{e}_i\\
\mathbf{x}_{l+1}=\mathbf{x}_l-\mathbf{y}_l\end{cases}.$$

We can easily show that $\mathbf{x}_1=\underset{i\in[n-1]}{\sum}\mathbf{y}_i$
and there are faces $F_1,F_2,\ldots,F_{n-1}$ of $\alpha_k^{d-1}$
such that
$$\begin{cases}
F_i\supseteq F_{i+1}\\
\mathbf{v}_{F_i}=\mathbf{x}_i\end{cases}.$$

\section{Decomposition theorems for polytope numbers\label{section4}} 
We have defined pointed triangulation in Section \ref{section2} and the geometric
description of polytope numbers in Section \ref{section3}. By
combining these concepts, we develop several methods to represent
polytope numbers by sums of simplex numbers, called decomposition
theorems. We also investigate relations between these decomposition
theorems.

\subsection{Decomposition theorem 1} Decomposition theorem 1 shows
that every $d$-polytope numbers can be decomposed into $d$-simplex
numbers.

\begin{theorem}[Decomposition theorem 1]\label{decompo1}
Let $P$ be a d-polytope with the $V(P)$-pointed triangulation
$\mathcal{C}_P$. Suppose that $P(n)$ are
the $V(P)$-polytope numbers. Then there are nonnegative integers
$a_1,a_2,\ldots,a_{d-1}$ such that
$$P(n)=\alpha^d(n)+\sum_{i\in[d-1]}a_i
\alpha^d(n-i).$$ \end{theorem}
\begin{pf}
Let $P$ be a $d$-polytope with the $V(P)$-pointed triangulation
$\mathcal{C}_{P}$. Then the geometric
description of polytope numbers furnishes $P(n)=\big|S_{P}(n)\big|$
where
$$S_{P}(n)=\bigcup_{\alpha^d\in\mathcal{F}_d(\mathcal{C}_{P})}S_{\alpha^d}(n).$$
Let
$\mathcal{F}_d(\mathcal{C}_{P})=\big\{\alpha^d_1,\alpha^d_2,\ldots,\alpha^d_s\big\}$.
Using Lemma \ref{decshell}, we can assign a shelling to $
\mathcal{F}_d(\mathcal{C}_{P})$ that satisfies the conditions of
Lemma \ref{decshell}. We assume that the indexes of elements in
$\mathcal{F}_d(\mathcal{C}_{P})$ equal those in the proof of Lemma
\ref{decshell}.

Let
$$\begin{cases}
\mathcal{P}_1=\big\{S_{\alpha^d_1}(n)\big\}\\
\mathcal{Q}_1=S_{\alpha^d_1}(n)\end{cases}.$$ Then
$$|\mathcal{Q}_1|=\alpha^d(n).$$
For $k\ge 2$ defining
$$\begin{cases}
\mathcal{P}_k=\mathcal{P}_{k-1}\cup
\big\{S_{\alpha^d_k}(n)\big\}\\
\mathcal{Q}_k=\mathcal{Q}_{k-1}\cup S_{\alpha^d_k}(n)\end{cases},$$
we inductively suppose that the number of facets of
$conv\big(S_{\alpha^d_k}(n)\big)$ in
$$\bigcup_{\alpha^d\in\mathcal{P}_{k-1}}conv\big(S_{\alpha^d}(n)\big),$$
denoted by $l_k$, satisfies $\l_k\in[d-1]$. Thus we need to
eliminate every point on those $l_k$ facets from
$conv\big(S_{\alpha^d_k}(n)\big)$ to compute $|\mathcal{Q}_k|$. For
this elimination, we apply successive $l_k$ facet-cuts to
$S_{\alpha^d_k}(n)$. Then
$$|\mathcal{Q}_k|=\alpha^d(n)+\sum_{i\in[k-1]}\alpha^d(n-l_{i+1}).$$
Continuing this operation until $k=m$ yields
$$|\mathcal{Q}_m|=\alpha^d(n)+\sum_{i\in[s-1]}\alpha^d(n-l_{i+1}).$$
Since $l_k\in[d-1]$ when $k\in[s]\setminus\{1\}$, there are nonnegative
integers $a_1,a_2,\ldots,a_{d-1}$ such that
$$\alpha^d(n)+\sum_{i\in[s-1]}\alpha^d(n-l_{i+1})=\alpha^d(n)+\sum_{i\in[d-1]}a_i\alpha^d(n-i).$$
By $\mathcal{Q}_m=S_{P }(n)$,
$$P(n)=\alpha^d(n)+\sum_{i\in[d-1]}a_i\alpha^d(n-i).$$\end{pf}

\begin{remark}For a $d$-polytope $P$ with the $V(P)$-pointed triangulation
$\mathcal{C}_P$, the number $P(n)$ is a polynomial in $n$ and it can
be determined by $d$ different values of $n$. Thus, whichever
shelling of $\mathcal{C}_P$ we may choose, the pointed triangulation
$\mathcal{C}_P$ uniquely determines the coefficients
$a_1,a_2,\ldots,a_{d-1}$. This means that a polynomial identity of
polytope numbers explains a geometric property of polytopes.
\end{remark}

\subsection{Decomposition theorem 2} Decomposition theorem 2 shows
that we can describe polytope numbers for a polytope by the facial
information of a pointed triangulation.

\begin{theorem}[Decomposition theorem 2]\label{decomp2}
Let $\mathcal{C}_{P}$ be the $V(P)$-pointed triangulation.
If $P(n)$ are the $V(P)$-polyopte numbers,
then
$$P(n)=\sum_{i\in[d-1]}(-1)^{d-i}b_i\alpha^i(n)$$
where $b_i$ is the number of $i$-simplexes $\alpha^i$ in
$\mathcal{C}_{P}$ such that
$$\begin{cases}
\mathbf{v}_{P}\in \alpha^i\\
\alpha^i\cap int(P)\neq\emptyset\end{cases}.$$
\end{theorem}

\begin{pf}
Let $\mathcal{C}_{P}$ be the $V(P)$-pointed triangulation.
We assume that the order on $d$-simplexes in
$\mathcal{C}_P$ is the same as that in Decomposition theorem 1. We can
naturally endow the $V(P_k)$-pointed triangulation of
$$P_k=\bigcup_{i\in[k]}conv\big(S_{\alpha^d_i}\big)$$ from $\mathcal{C}_{P}$ for $k\in[s]$.
Even when $P_k$ is not convex, this triangulation is possible by the
definition of pointed triangulation. We denote by
$\mathcal{C}_{P_k}$ the endowed $V(P_k)$-pointed triangulation of $P_k$.

Let $S_{P_k}(n)=\underset{i\in[k]}{\bigcup} S_{\alpha^d_i}(n)$. If $k=1$ then
$$\big|S_{P_1}(n)\big|=\alpha^d(n).$$ Inductively, we assume that
$$\big|S_{P_k}(n)\big|=\sum_{i\in[d]}(-1)^{d-i}b_{k,i}\alpha^d_i(n),$$
where $b_{k,i}$ is the number of $\alpha^i\in\mathcal{C}_{P_k}$ such
that
$$\begin{cases}
\mathbf{v}_{P}\in\alpha^i\\
\alpha^i\cap
int\bigg(\underset{j\in[k]}{\bigcup}\alpha^d_j\bigg)\neq\emptyset\end{cases}.$$
We claim that 
\begin{equation}\label{decomposition2}
\big|S_{P_{k+1}}(n)\big|=\sum_{i\in[d]}(-1)^{d-i}b_{k+1,i}\alpha^d_i(n).
\end{equation}
If either $d=0$ or $d=1$, then the identity (\ref{decomposition2}) is obvious.
Suppose that $d\ge 2$. For $$F^{d-1}_{k+1}=\alpha^d_{k+1}\cap
\big|\mathcal{C}_{P_k}\big|$$ let $\mathcal{C}_{F^{d-1}_{k+1}}$ be
the endowed $V\big(F^{d-1}_{k+1}\big)$-pointed triangulation of $F^{d-1}_{k+1}$ from
$\mathcal{C}_{P_k}$. By the definition of pointed triangulation,
there are $(d-2)$-simplexes
$$\alpha^{d-2}_{k+1,1},\alpha^{d-2}_{k+1,2},\ldots,\alpha^{d-2}_{k+1,l_{k+1}}$$
such that
$$F^{d-1}_{k+1}=\bigcup_{i\in[l_{k+1}]}\alpha^{d-2}_{k+1,i}.$$
Thus $\mathcal{C}_{F^{d-1}_{k+1}}$ is the $V\big(F^{d-1}_{k+1}\big)$-pointed triangulation of
$F_{k+1}^{d-1}$. By induction on
$d$
$$\bigg|\bigcup_{i\in[l_{k+1}]}S_{\alpha^{d-1}_{k+1,i}}(n)\bigg|
=\sum_{i\in[d-1]}(-1)^{d-i}c_{k+1,i}\alpha^i(n),$$ where $c_{k+1,i}$
is the number of $i$-simplexes $\alpha^i$ in
$\mathcal{C}_{F^{d-1}_{k+1}}$ such that
$$\begin{cases}
\mathbf{v}_{P}\in\alpha^i\\
\alpha^i\cap int \big(F^{d-1}_{k+1}
\big)\neq\emptyset\end{cases}.$$ Hence
\begin{align*}
\big|S_{P_{k+1}}(n)\big|&
=\big|S_{P_k}(n)\big|+\big|\alpha^d_{k+1}(n)\big|
-\bigg|\bigcup_{i\in[l_{k+1}]}\alpha^{d-1}_{k+1,i}(n)\bigg|\\
                   &=\sum_{i\in[d]}(-1)^{d-i}
                   b_{k,i}\alpha^i(n)+\alpha^d(n)-\sum_{i\in[d-1]}(-1)^{d-1-i}c_{k+1,i}\alpha^i(n)\\
                   &=(b_{k,d}+1)\alpha^d(n)+\sum_{i\in[d-1]}(-1)^{d-i}(b_{k,\,i}+c_{k+1,i})\alpha^i(n).
\end{align*}
Since $b_{k,d}+1$ is the number of $d$-simplexes $\alpha^d$ in
$\mathcal{C}_{P_{k+1}}$ such that
$$\begin{cases}
\mathbf{v}_{P}\in \alpha^d\\
\alpha^d\cap
int\bigg(\underset{j\in[k+1]}{\bigcup}\alpha^{d}_{k+1}\bigg)\neq\emptyset\end{cases}$$
and for $i\in[d-1]$ the number $b_{k,i}+c_{k+1,i}$ is the
number of $i$-simplexes $\alpha^i$ in $\mathcal{C}_{k+1}$ such that
$$\begin{cases}
\mathbf{v}_{P}\in\alpha^i\\
\alpha^i\cap
int\bigg(\underset{i\in[k+1]}{\bigcup}\alpha^{d}_{k+1}\bigg)\neq\emptyset\end{cases},$$
we prove our assertion.

If we continue this process until $k=s$, then
$$P(n)=\sum_{i\in[d]}(-1)^{d-i}b_i\alpha^i(n)$$
where $b_i$ is the number of $\alpha^i$ in $P$ such that
$$\begin{cases}
\mathbf{v}_{P}\in \alpha^i\\
\alpha^i\cap int(P)\neq\emptyset\end{cases}.$$
\end{pf}

\subsection{Decomposition theorem 3} The main idea of Decomposition
theorem 3 is to represent a pointed triangulation by a disjoint
union and then to apply this representation to polytope numbers.

\begin{theorem}[Decomposition theorem 3-1]\label{dec3-1}
Let $\mathcal{C}_{P}$ be the $V(P)$-pointed triangulation
and
$$\mathcal{F}_k(\mathcal{C}_{P},\mathbf{v}_{P})=\big\{\alpha^k\,\big|\,\mathbf{v}_{P}\in
\alpha^k, \alpha^k\in\mathcal{F}_k (\mathcal{C}_{P})\big\}.$$ Then the $V(P)$-polytope numbers are represented by
$$P(n)=\sum_{k\in[d]}c_k\alpha^k(n-k)$$ where
$c_k=\big|\mathcal{F}_k(\mathcal{C}_{P},\mathbf{v}_{P})\big|$.
\end{theorem}
\begin{pf}
For each $\alpha^k\in\mathcal{F}_k(\mathcal{C}_{P},\mathbf{v}_{P})$
there are exactly $k$ facets $\alpha^{k-1}$ of $\alpha^k$ that are in
$\mathcal{F}_{k-1}(\mathcal{C}_{P},\mathbf{v}_{P})$ such that $\mathbf{v}_{P}\in \alpha^{k-1}$.
Correlating $\alpha^k(n)$ to $S_{\alpha^k}(n)$, we apply successive
$k$ facet-cuts to $S_{\alpha^k}(n)$. Then these successive facet-cuts
change $S_{\alpha^k}(n)$ into $S^k_{\alpha^k}(n)$ where
$$\big|S^k_{\alpha^k}(n)\big|=\alpha^k(n-k).$$ Each pair of distinct faces $\alpha^i$ and $\alpha^j$ of $\mathcal{C}_P$ contained in
$\mathcal{F}_i(\mathcal{C}_{P},\mathbf{v}_{P})$ satisfies
$$S^i_{\alpha^i}(n)\cap S^j_{\alpha^j}(n)=\emptyset,$$ thus
\begin{align*}
P(n)&=\sum_{k\in[d]}\sum_{\alpha^k\in\mathcal{F}_k(\mathcal{C}_{P},\mathbf{v}_{P})}
 S^k_{\alpha^k}(n)\\
 &=\sum_{k\in[d]}c_k\alpha^k(n-k).\end{align*}
\end{pf}

Instead of considering all of the faces in $\mathcal{C}_P$, if we
consider only those faces of $\mathcal{C}_{P}$ that have nonempty
intersection with $int(P)$, then we obtain the following corollary.
\begin{corollary}[Decomposition theorem 3-2]\label{dec3-2}
Let
$$\begin{cases}
\mathcal{F}_{-1}^1(\mathcal{C}_{F},\mathbf{v}_{P})=\emptyset\\
\mathcal{F}_{0}^1(\mathcal{C}_{F},\mathbf{v}_{P})=\{\mathbf{v}_P\}\end{cases},$$
and for $k\in[d]$ let
$$\mathcal{F}_{k}^1(\mathcal{C}_{F},\mathbf{v}_{P})=\Big\{\alpha^k\,
\Big|\,\alpha^k\in\mathcal{F}_k(\mathcal{C}_{P},\mathbf{v}_{P}),
\alpha^k\cap int(P)\neq\emptyset\Big\}.$$ For each
$\alpha^k\in\mathcal{F}_{k}^1(\mathcal{C}_{F},\mathbf{v}_{P})$
with $k\in[d]_0$ let $f(\alpha^k)$ be the number of facets of
$\alpha^k$ that are in
$\mathcal{F}_{k-1}^1(\mathcal{C}_{F},\mathbf{v}_{P})$. Then the $V(P)$-polytope numbers are represented by
$$P(n)=\sum_{k\in[d]_0}\sum_{\alpha^k\in\mathcal{F}_{k}^1(\mathcal{C}_{F},\mathbf{v}_{P})}
\alpha^k\big(n-f(\alpha^k)\big).$$
\end{corollary}

\subsection{Decomposition theorem 4}
Let $\mathcal{C}_P$ be a pointed triangulation of $P$. Then
$$P=\biguplus_{k\in[d]_0}\biguplus_{\alpha^k\in\mathcal{F}_k(\mathcal{C}_P)}relint(\alpha^k).$$
Using this decomposition, we consider Decomposition theorem $4$.
\begin{theorem}[Decomposition theorem 4]\label{dec4}
Let $P$ be a d-polytope with the $V(P)$-pointed triangulation
$\mathcal{C}_{P}$ and for $k\in[d]_0$
let $d_k$ be the number of $k$-simplexes in $\mathcal{C}_{P}$.
Then the 
$V(P)$-polytope numbers are represented by
$$P(n)=\sum_{k\in[d]}d_k\alpha^k\big(n-(k+1)\big).$$
\end{theorem}
\begin{pf}
Each pair of two distinct faces $\alpha^{k_1}$ and $\alpha^{k_2}$ of
$\mathcal{C}_{P}$ satisfies
$$relint(\alpha^{k_1})\cap
relint(\alpha^{k_2})=\emptyset,$$ thus
\begin{equation}\label{disjoint}
P=\biguplus_{k\in[d]_0}\biguplus_{\alpha^k\in\mathcal{C}_{P}}relint(\alpha^k).
\end{equation}
The identity (\ref{disjoint}) yields
$$S_P(n)=\biguplus_{k\in[d]_0}\biguplus_{\alpha^k\in\mathcal{F}(\mathcal{C}_P)}
S_{\alpha^k}(n)^{\sharp}.$$ Since the interior $k$-simplex numbers $\alpha^k(n)^{\sharp}$ satisfy
$$\alpha^k(n)^{\sharp}=\alpha^k\big(n-(k+1)\big),$$
if we denote $d_k=\big|\mathcal{F}_k(\mathcal{C}_P)\big|$ then
\begin{align*}
P(n)&=\big|S_P(n)\big|=\bigg|\biguplus_{k\in[d]_0}\biguplus_{\alpha^k\in\mathcal{F}(\mathcal{C}_P)}
S_{\alpha^k}(n)^{\sharp}\bigg|\\
&=\sum_{k\in[d]_0}\sum_{\alpha^k\in
\mathcal{F}_k(\mathcal{C}_P)}\alpha^k\big(n-(k+1)\big)\\
&=\sum_{k\in[d]_0}d_k \alpha^k\big(n-(k+1)\big).
\end{align*}
\end{pf}

\subsection{Relations between decomposition theorems} Decomposition
theorems provides several methods to decompose polytope numbers into
simplex numbers. Thus relations between decomposition theorems may
exist. We consider such relations in this subsection.

\subsubsection{Decomposition theorem 1 and other decomposition
theorems} Decomposition theorem 1 represents $d$-polytope numbers by
sums of $d$-simplex numbers, whereas other decomposition theorems
represent $d$-polytope numbers by sums of various dimensional
simplex numbers. Thus we use the relation
$\alpha^{d-1}(n)=\alpha^d(n)-\alpha^d(n-1)$ to derive Decomposition
theorem 1 from other decomposition theorems.

The equation $\alpha^{d-1}(n)=\alpha^d(n)-\alpha^d(n-1)$ produces
\begin{align}
&\alpha^i(n)=\sum_{j\in[d-i]_0}(-1)^j{d-i\choose j}\alpha^d(n-j),\label{aidentity1}\\
&\alpha^k(n-k)=\sum_{j\in[d-k]_0}(-1)^j{d-k\choose
j}\alpha^d(n-k-j),\label{aidentity2}\\
&\alpha^k\big(n-(k+1)\big)=\sum_{j\in[d-k]_0}(-1)^j{d-k\choose
j}\alpha^d\big(n-(k+1)-j\big)\label{aidentity3}.
\end{align}
Therefore for a $d$-polytope $P$ the identities
(\ref{aidentity1})--(\ref{aidentity3}) change decomposition forms of
$P(n)$ in Decomposition theorems 2, 3-1, and 4 into
\begin{align*}
P(n)&=\sum_{i\in[d]}(-1)^{d-i}b_i a^i(n)\\
&=\sum_{i\in[d]}(-1)^{d-i}b_i\bigg(\sum_{j\in[d-i]}(-1)^j{d-i\choose
j}\alpha^d(n-j)\bigg)\\
&=\sum_{j\in[d-1]_0}\bigg(\sum_{i\in[d-j]}(-1)^{d-i-j}b_i{d-i\choose
j}\bigg)\alpha^d(n-j),
\end{align*}
\begin{align*}
P(n)&=\sum_{k\in[d]_0} c_k\alpha^k(n-k)\\
      &=\sum_{k\in[d]_0} c_k\bigg(\sum_{j\in[d-k]_0}(-1)^j{d-k\choose
      j}\alpha^d(n-k-j)\bigg) \\
      &=\sum_{l\in[d]_0}\bigg(\sum_{k\in[l]_0}c_k(-1)^{l-k}{d-k\choose
      l-k}\bigg)\alpha^d(n-l),
\end{align*}
and
\begin{align*}
P(n)&=\sum_{k\in[d]_0} d_k\alpha^k(n-(k+1))\\
      &=\sum_{k\in[d]_0} d_k\bigg(\sum_{j\in[d-k]_0}(-1)^j{d-k\choose
      j}\alpha^d\big(n-(k+1)-j\big)\bigg)\ \\
      &=\sum_{l\in[d+1]}\bigg(\sum_{k\in[l-1]_0}d_k(-1)^{l-1-k}{d-k\choose
      l-1-k}\bigg)\alpha^d(n-l),
\end{align*}
respectively.

\subsubsection{Decomposition theorems 2 and 3-2}
For a polytopal complex $\mathcal{C}_P$, Decomposition theorems 2
and 3-2 use the same facial information of $\mathcal{C}_P$. Hence we can expect a
relation between these two decomposition theorems. We verify such a relation here.

\medskip
For a finite poset $\mathcal{P}$ we define the \emph{zeta function}
$\zeta$ of $\mathcal{P}$ to be
$$\zeta(x,y)=1\ \text{for each}\ \{x,y\}\subseteq\mathcal{P}\ \text{with}\ x\leq y.$$
Let the inverse function of $\zeta$ be $\mu$ called the
\emph{M\"{o}bius function} of $\mathcal{P}$. Then
$$
\mu(x,y)=\begin{cases}
1&\text{when}\ x=y\\
-\underset{x\leq z<y}{\sum}\mu(x,z)&\text{when}\ x<y\end{cases}.
$$
\begin{theorem}[M\"{o}bius inversion formula \cite{stanley}]\label{mobius}
Let $\mathcal{P}$ be a finite poset and $f:\mathcal{P}\rightarrow
\mathbb{C}$ and $g:\mathcal{P}\rightarrow
\mathbb{C}$. Then
$$g(x)=\sum_{y\leq x}f(y)\ \text{when}\ x\in \mathcal{P}$$
if and only if
$$f(x)=\sum_{y\leq x}g(y)\mu(y,x)\ \text{when}\ x\in \mathcal{P}.$$
\end{theorem}

\begin{theorem}[\cite{stanley}]\label{proposition4}
Let $\mathcal{C}$ be a polytopal complex and
$\mathcal{P}=\mathcal{P}(\mathcal{C})$ be the poset on
$\big(\mathcal{C}\setminus\{\emptyset\}\big)
\cup\big\{|\mathcal{C}|\big\}$, ordered by $$F_i\leq F_j\ \text{if}\
F_i\subseteq F_j.$$ Then
$$
\mu_{\mathcal{P}}(F_1,F_2)=
\begin{cases}
0,&\text{if $F_2=|\mathcal{C}|$ and $F_1$ lies on
the}\\
&\text{boundary of $|\mathcal{C}|$}\\
(-1)^{rank(F_2)-rank(F_1)},&\text{otherwise}
\end{cases}.$$
\end{theorem}

Let $P$ be a $d$-polytope with the $V(P)$-pointed triangulation
$\mathcal{C}_{P}$ and $\mathcal{C}_{P^{d-1}}$
be the polytopal complex formed by the faces of $\mathcal{C}_P$ that
do not contain $\mathbf{v}_{P}$. We define
\begin{equation}\begin{cases}\label{complex23}
{\mathcal{C}_{P^{d-1}}}(n)^{\sharp}=0\\
\mathcal{C}_{P^{d-1}}(n)={\mathcal{C}_{P^{d-1}}}(n)^{\sharp}
+\underset{k\in[d-1]_0}{\sum}\underset{\alpha^k\in\mathcal{C}_{P^{d-1}}}{\sum}
{\alpha^k(n)}^{\sharp}
\end{cases}.\end{equation}
Note that the identity (\ref{complex23}) is a restatement of the
decomposition form of Decomposition theorem 3-2. If we let
$\mathcal{P}=\mathcal{P}(\mathcal{C}_{P^{d-1}})$, then Theorems
\ref{mobius} and \ref{proposition4} supply
$${\mathcal{C}_{P^{d-1}}(n)}^{\sharp}=\mathcal{C}_{P^{d-1}}(n)
+\sum_{k\in[d-1]_0}(-1)^{d-k} b_{k+1}\alpha^k(n)$$ where $b_{k+1}$ is
the number of $(k+1)$-simplexes that contain $\mathbf{v}_P$ and are
not on the boundary of $P$. Since
${\mathcal{C}_{P^{d-1}}(n)}^{\sharp}=0$, we obtain
$$\mathcal{C}_{P^{d-1}}(n)=\sum_{k\in[d-1]_0}(-1)^{d-1-k}b_{k+1}\alpha^k(n).$$
Moreover,
$$\begin{cases}
P (n)=\underset{i\in[n]}{\sum}\mathcal{C}_{P^{d-1}}(i)\\
\alpha^{k+1}(n)=\underset{i\in[n]}{\sum}\alpha^k(i)\end{cases}$$ by definition,
thus
$$P(n)=\underset{k\in[d]}{\sum}(-1)^{d-k}b_k\alpha^k(n),$$ which is the decomposition form
of Decomposition theorem 2.

Instead of triangulation, we consider the polytopal complex
$\mathcal{C}_{Q^{d-1}}$ formed by the faces of $P$ that do not
contain $\mathbf{v}_{P}$. We define
$$
\begin{cases}\mathcal{C}_{Q^{d-1}}(n)^{\sharp}=0\\
\mathcal{C}_{Q^{d-1}}(n)=\mathcal{C}_{Q^{d-1}}(n)^{\sharp}+
\underset{k\in[d-1]_0}{\sum}\underset{F^k\in\mathcal{C}_{Q^{d-1}}}{\sum}F^k(n)^{\sharp}
\end{cases}.$$ The method used in the triangulation case yields
$$\mathcal{C}_{Q^{d-1}}(n)=\sum_{k\in[d-1]}(-1)^{d-1-k}
\sum_{F^k\in\mathcal{C}_{Q^{d-1}}}F^k(n).$$ Since $P
(n)=\underset{i\in[n]}{\sum}\mathcal{C}_{Q^{d-1}}(i)$, if we define
$F^{k,s}(n)=\underset{i\in[n]}{\sum}F^k(i)$ then
$$P(n)=\sum_{k\in[d]}(-1)^{d-k}
\sum_{F^k\in\mathcal{C}_{Q^{d-1}}}F^{k,\,s}(n).$$

\subsection{Computations of coefficients in decomposition theorems}
Let $P$ be a $d$-polytope. By the definition of polytope numbers,
for each $n$ the number $P(n)$ is a polynomial of $n$. Therefore a
finite number of values of $P(n)$ allow us to compute the
coefficients in decomposition forms of $P(n)$ in decomposition theorems.
We perform such computations.

\medskip
By Decomposition theorem 1
$$P(n)=\sum_{i\in[d-1]}a_i\alpha^d(n-i),$$
thus
$$\begin{pmatrix}
P(1)\\P(2)\\ \cdots\\ P(d)
\end{pmatrix}
=\begin{pmatrix} \alpha^d(1)&0&\cdots&0\\
\alpha^d(2)&\alpha^d(1)&\cdots&0\\
\cdots&\cdots&\cdots&\cdots\\
\alpha^d(d)&\alpha^d(d-1)&\cdots&\alpha^d(1)
\end{pmatrix}\begin{pmatrix}
a_{0}\\a_{1}\\ \cdots\\a_{d-1}
\end{pmatrix}.
$$
If we let $$A=\begin{pmatrix} \alpha^d(1)&0&\cdots&0\\
\alpha^d(2)&\alpha^d(1)&\cdots&0\\
\cdots&\cdots&\cdots&\cdots\\
\alpha^d(d)&\alpha^d(d-1)&\cdots&\alpha^d(1)
\end{pmatrix},$$
then the Gaussian elimination produces
$$A^{-1}=\begin{pmatrix}
(-1)^0{d+1\choose 0}&0&\cdots&0\\
(-1)^1{d+1\choose 1}&(-1)^0{d+1\choose 0}&\cdots&0\\
\cdots&\cdots&\cdots&\cdots\\
(-1)^{d-1}{d+1\choose d-1}&(-1)^{d-2}{d+1\choose
d-2}&\cdots&(-1)^0{d+1\choose 0}
\end{pmatrix}.$$
Therefore
$$\begin{pmatrix}
a_{0}\\a_{1}\\ \cdots\\a_{d-1}
\end{pmatrix}
=
\begin{pmatrix}
(-1)^0{d+1\choose 0}&0&\cdots&0\\
(-1)^1{d+1\choose 1}&(-1)^0{d+1\choose 0}&\cdots&0\\
\cdots&\cdots&\cdots&\cdots\\
(-1)^{d-1}{d+1\choose d-1}&(-1)^{d-2}{d+1\choose
d-2}&\cdots&(-1)^0{d+1\choose 0}
\end{pmatrix}
\begin{pmatrix}
P(1)\\P(2)\\ \cdots\\ P(d)
\end{pmatrix}.
$$

By Decomposition theorem 2
$$P(n)=\sum_{i\in[d]}(-1)^{d-i}b_i\alpha^i(n),$$
which gives
$$\begin{pmatrix}
P(1)\\P(2)\\ \cdots\\P(d)
\end{pmatrix}=
\begin{pmatrix}
\alpha^1(1)&\alpha^2(1)&\cdots&\alpha^d(1)\\
\alpha^1(2)&\alpha^2(2)&\cdots&\alpha^d(2)\\
\cdots&\cdots&\cdots&\cdots\\
\alpha^1(d)&\alpha^2(d)&\cdots&\alpha^d(d)
\end{pmatrix}
\begin{pmatrix}
(-1)^{d-1}b_1\\(-1)^{d-2}b_2\\ \cdots\\(-1)^0 b_d
\end{pmatrix}.$$
If we denote
$$B=
\begin{pmatrix}
\alpha^1(1)&\alpha^2(1)&\cdots&\alpha^d(1)\\
\alpha^1(2)&\alpha^2(2)&\cdots&\alpha^d(2)\\
\cdots&\cdots&\cdots&\cdots\\
\alpha^1(d)&\alpha^2(d)&\cdots&\alpha^d(d)
\end{pmatrix}$$
then
$$\begin{pmatrix} (-1)^{d-1}b_1\\(-1)^{d-2}b_2\\
\cdots\\(-1)^1b_{d-1}\\(-1)^0 b_d
\end{pmatrix}\\
=B^{-1}\begin{pmatrix} P(1)\\P(2)\\ \cdots\\P(d-1)\\P(d)
\end{pmatrix}$$ where the matrix $B^{-1}$ is
$$\begin{pmatrix}
(-1)^0({1\choose 0}+{2\choose 1})&(-1)^{-1}{2\choose 0}&0&\cdots&0\\
(-1)^1({2\choose 1}+{3\choose 2})&(-1)^0({2\choose 0}+{3\choose
1})&(-1)^{-1}{3\choose 0}&\cdots&0\\
\cdots&\cdots&\cdots&\cdots&\cdots\\
(-1)^{d-2}({d-1\choose d-2}+{d\choose
d-1})&\cdots&\cdots&(-1)^0({d-1\choose 0}+{d\choose
1})&(-1)^{-1}{d\choose 0}\\
(-1)^{d-1}{d\choose d-1}&\cdots&\cdots&(-1)^1{d\choose
1}&(-1)^0{d\choose 0}
\end{pmatrix}.$$

Decomposition theorem 3 yields
$$P(n)=\sum_{i\in[d]}c_i\alpha^d(n-i).$$
Similarly, if we write
$$\begin{pmatrix}
P(1)\\P(2)\\ \cdots\\P(d+1)
\end{pmatrix}=\begin{pmatrix}
\alpha^0(1)&0&\cdots&0\\
\alpha^0(2)&\alpha^1(1)&\cdots&0\\
\cdots&\cdots&\cdots&\cdots\\
\alpha^0(d+1)&\alpha^1(d)&\cdots&\alpha^d(1)
\end{pmatrix}
\begin{pmatrix}
c_0\\ c_1\\ \cdots\\c_d
\end{pmatrix}$$
then
$$
\begin{pmatrix}
c_0\\ c_1\\ \cdots\\c_d
\end{pmatrix}
=\begin{pmatrix}
(-1)^0{0\choose 0}&0&\cdots&0\\
(-1)^1{1\choose 1}&(-1)^0{1\choose 0}&\cdots&0\\
\cdots&\cdots&\cdots&\cdots\\
(-1)^d{d\choose d}&(-1)^{d-1}{d\choose d-1}&\cdots&(-1)^0{d\choose
0}
\end{pmatrix}\begin{pmatrix} P(1)\\P(2)\\ \cdots\\P(d+1)
\end{pmatrix}.
$$

\section{Illustrations of decomposition theorems\label{section5}}
We apply decomposition theorems to both regular polytopes and the
product of simplexes. This application gives new interpretations of
some known combinatorial identities and derives new combinatorial
identities.

\subsection{Decomposition theorem 1}

\subsubsection{Cross polytope $\beta^d$} For fixed numbers $a_i\in\{1,-1\}$ let
$\alpha^d(a_1,a_2,\ldots,a_{d-1})$ be the $d$-simplex with the
vertex set
$$\{\mathbf{e}_d,-\mathbf{e}_d,a_1\mathbf{e}_1,a_2\mathbf{e}_2,\ldots,
a_{d-1}\mathbf{e}_{d-1}\},$$ and let
$\mathcal{S}\big(\beta_s^d\big)$ be the set of such $d$-simplexes.
Then, by a simple reasoning,
$$\beta^d=\bigcup_{\alpha^d\in\mathcal{S}(\beta^d)}\alpha^d.$$
Moreover, the $d$-simplexes in $\mathcal{S}\big(\beta_s^d\big)$ form
a pointed triangulation $\mathcal{C}_{\beta_s^d}$ such that
$$\begin{cases}\big|\mathcal{C}_{\beta_s^d}\big|=\beta^d\\
\mathcal{C}_{\beta_s^d}
=\underset{\alpha^d\in\mathcal{S}(\beta^d)}{\bigoplus}\alpha^d\end{cases}.$$
We call $\mathcal{C}_{\beta_s^d}$ the \emph{standard pointed
triangulation} of $\beta_s^d$ and assume that $\beta_s^d(n)$ are
formed by $\mathcal{C}_{\beta_s^d}$.

For $i\in[d-1]_0$ let $\mathcal{S}_i(\beta_s^d)$ be the set of
$\alpha^d(a_1,a_2,\ldots,a_{d-1})$ in $\mathcal{S}(\beta_s^d)$ such
that the number of $a_j$ satisfying $a_j=-1$ is $i$. Then for each
$i$ the set $\mathcal{S}_i(\beta_s^d)$ satisfies
$$\begin{cases}
\big|\mathcal{S}_i(\beta_s^d)\big|={d-1\choose i}\\
\mathcal{S}(\beta_s^d)=\underset{i\in[d-1]}{\biguplus}\mathcal{S}_i(\beta_s^d)\end{cases}.$$
We say that a $d$-simplex
$\alpha^d\big(a_1^1,a_2^1,\ldots,a_{d-1}^1\big)$ in
$\mathcal{S}(\beta_s^d)$ is \emph{adjacent} to another $d$-simplex
$\alpha^d\big(a_1^2,a_2^2,\ldots,a_{d-1}^2\big)$ in
$\mathcal{S}(\beta_s^d)$ if
$$\big(a_1^2,a_2^2,\ldots,a_{d-1}^2\big)-\big(a_1^1,a_2^1,\ldots,a_{d-1}^1\big)=2\mathbf{e}_j$$
for some $j\in[d-1]$. By the definition of adjacency
$$dim\Big(\alpha^d\big(a_1^1,a_2^1,\ldots,\,a_{d-1}^1\big)\cap
\alpha^d\big(a_1^2,a_2^2,\ldots,a_{d-1}^2\big)\Big)=d-1,$$ and if
$\alpha^d\big(a_1^1,a_2^1,\ldots,a_{d-1}^1\big)\in\mathcal{S}_i(\beta_s^d)$
then
$\alpha^d\big(a_1^2,a_2^2,\ldots,a_{d-1}^2\big)\in\mathcal{S}_{i+1}(\beta_s^d)$.
In addition, if
$\alpha^d(a_1,a_2,\ldots,a_{d_1})\in\mathcal{S}_i(\beta_s^d)$ then
the number of $d$-simplexes in $\mathcal{S}(\beta_s^d)$ that are
adjacent to $\alpha^d(a_1,a_2,\ldots,a_{d-1})$ is $i$.

For $\alpha^d\in\mathcal{S}(\beta_s^d)$ let $A(\alpha^d)$ be the set
of $d$-simplexes in $\mathcal{S}(\beta_s^d)$ that are adjacent to
$\alpha^d$ and
$$C_{\alpha^d}(n)=S_{\alpha^d}(n)\setminus\bigcup_{\alpha_1^d\in
A(\alpha^d)}S_{\alpha_1^d}(n).$$ Simply speaking, $C_{\alpha^d}(n)$
is formed by successive facet-cuts on $S_{\alpha^d}(n)$. By the
geometric description of polytope numbers
$$S_{\beta_s^d}(n)=\biguplus_{\alpha_d\in
\mathcal{S}(\beta_s^d)}C_{\alpha^d}(n),$$ and if
$\alpha^d\in\mathcal{S}_i$ then successive $i$ facet-cuts on
$S_{\alpha^d}$ yield
$$\big|C_{\alpha^d}(n)\big|=\alpha^d(n-i).$$
Therefore
\begin{align*}
\beta^d(n)&=\beta_s^d(n)=\big|S_{\beta_s^d}(n)\big|\\
&=\bigg|\biguplus_{\alpha_d\in
\mathcal{S}(\beta_s^d)}C_{\alpha^d}(n)\bigg|
=\bigg|\biguplus_{i\in[d-1]_0}\biguplus_{\alpha^d\in\mathcal{S}_i(\beta_s^d)}C_{\alpha^d}(n)\bigg|\\
&=\sum_{i\in[d-1]_0}\sum_{\alpha^d\in\mathcal{S}_i(\beta_s^d)}\big|C_{\alpha^d}(n)\big|
=\sum_{i\in[d-1]_0}\sum_{\alpha^d\in\mathcal{S}_i(\beta_s^d)}\alpha^d(n-i)\\
&=\sum_{i\in[d-1]_0}{d-1\choose i}\alpha^d(n-i).
\end{align*}
This computation of $\beta^d(n)$ coincides with Kim's computation of
$\beta^d(n)$ \cite{kim}.

\subsubsection{Measure polytope $\gamma^d$}

Let $S_d$ be the set of permutations on the set $[d]$.
We denote by $[a_1,a_2,\ldots,a_n]$ a permutation on
$[d]$. For the permutation $[a_1,a_2,\ldots,a_d]$ we
define $\alpha^d[ a_1,a_2,\ldots,a_d]$ to be the simplex with the
vertex set
$$\{\mathbf{0},\mathbf{e}_{a_1},\mathbf{e}_{a_1}+\mathbf{e}_{a_2},\ldots,
\mathbf{e}_{a_1}+\mathbf{e}_{a_2}+\cdots+\mathbf{e}_{a_d}=\mathbf{1}
\}.$$ If we define
$$\mathcal{S}(\gamma_s^d)=\big\{\alpha^d[a_1,a_2,\ldots,a_d]\,\big|\,[a_1,a_2,\ldots,a_d]\in S_d\big\},$$
then, by a simple reasoning,
$$\gamma_s^d=\bigcup_{\alpha^d\in\mathcal{S}(\gamma_s^d)}\alpha^d.$$ Stanley gave this
representation of $\gamma_s^d$ \cite{stanley}. We can easily show
that the $d$-simplexes in $\mathcal{S}(\gamma_s^d)$ form a pointed
triangulation $\mathcal{C}_{\gamma_s^d}$ such that
$$\begin{cases}
|\mathcal{C}_{\gamma_s^d}|=\gamma_s^d\\
\mathcal{C}_{\gamma^d}=\underset{\alpha^d\in\mathcal{S}(\gamma_s^d)}{\bigoplus}\alpha^d\end{cases}.$$
We call $\mathcal{C}_{\gamma_s^d}$ the \emph{standard pointed
triangulation} of $\gamma_s^d$ and assume that $\gamma^d(n)$ are
formed by this pointed triangulation

For $i\in[d-1]_0$ let $\mathcal{S}_i(\gamma_s^d)$ be the set of
$\alpha^d[ a_1,a_2,\dots,a_d]$ in $\mathcal{S}(\gamma_s^d)$ such
that for $j\in[d-1]$ the number of $j$ obeying $a_j>a_{j+1}$
is $i$. Each of the sets $\mathcal{S}_i(\gamma_s^d)$ satisfies
$$\begin{cases}\big|\mathcal{S}_i(\gamma_s^d)\big|=\bigg\langle{d\atop
i}\bigg\rangle\\
\mathcal{S}(\gamma_s^d)=\underset{i\in[d-1]_0}{\bigoplus}\mathcal{S}_i(\gamma_s^d)\end{cases}$$
where $\Big\langle{d\atop i}\Big\rangle$ is the Eulerian number. We
say that an $\alpha^d\big[a^1_1,a^1_2,\ldots,a^1_d\big]$ is
\emph{adjacent} to another $\alpha^d\big[
a^2_1,a^2_2,\ldots,a^2_d\big]$ if
$$\big[a_1^2,a_2^2,\ldots,a_d^2\big]
=\big[a_1^1,\ldots,a_{j-1}^1,a_{j+1}^1,a_{j}^1,a_{j+2}^1,\ldots,a_d^1\big]$$
where $a_j^1<a_{j+1}^1$. By the definition of adjacency
$$dim\Big(\alpha^d\big[a_1^1,a_2^1,\ldots,a_{d}^1\big]\bigcap
\alpha^d\big[ a_1^2,a_2^2,\ldots,a_{d}^2\big]\Big)=d-1,$$ and if
$\alpha^d[ a_1^1,a_2^1,\ldots,a_{d}^1]\in\mathcal{S}_i(\gamma_s^d)$
then $\alpha^d[
a_1^2,a_2^2,\ldots,a_{d}^2]\in\mathcal{S}_{i+1}(\gamma_s^d)$. In
addition, if
$\alpha^d[a_1,a_2,\ldots,a_{d}]\in\mathcal{S}_i(\gamma_s^d)$ then
the number of $d$-simplexes in $\mathcal{S}(\gamma_s^d)$ that are
adjacent to $\alpha^d[ a_1,a_2,\ldots,a_{d}]$ is $i$.

For $\alpha^d\in\mathcal{S}(\gamma_s^d)$ let $\mathcal{A}(\alpha^d)$
be the set of $d$-simplexes in $\mathcal{S}(\gamma_s^d)$ that are
adjacent to $\alpha^d$ and
$$C_{\alpha^d}(n)=S_{\alpha^d}(n)\setminus\bigcup_{\alpha_1^d\in
\mathcal{A}(\alpha^d)}S_{\alpha_1^d}(n).$$ Simply speaking,
$C_{\alpha^d}(n)$ is formed by successive facet-cuts on
$S_{\alpha^d}(n)$. By the geometric description of polytope numbers
$$S_{\gamma_s^d}(n)=\biguplus_{\alpha_d\in
\mathcal{S}(\gamma_s^d)}C_{\alpha^d}(n),$$ and if
$\alpha^d\in\mathcal{S}_i(\gamma_s^d)$ then successive $i$ facet-cuts on $S_{\alpha^d}(n)$ yield
$$\big|C_{\alpha^d}(n)\big|=\alpha^d(n-i).$$
Therefore
\begin{align}
\gamma^d(n)&=\gamma_s^d=\big|S_{\gamma_s^d}(n)\big|\label{dmeasure}\\
&=\Bigg|\biguplus_{\alpha_d\in
\mathcal{S}(\gamma_s^d)}C_{\alpha^d}(n)\Bigg|
=\bigg|\biguplus_{i\in[d-1]_0}\biguplus_{\alpha^d\in\mathcal{S}_i(\gamma_s^d)}C_{\alpha^d}(n)\bigg|\nonumber\\
&=\sum_{i\in[d-1]_0}\sum_{\alpha^d\in\mathcal{S}_i(\gamma_s^d)}\big|C_{\alpha^d}(n)\big|
=\sum_{i\in[d-1]_0}\sum_{\alpha^d\in\mathcal{S}_i(\gamma_s^d)}\alpha^d(n-i)\nonumber\\
&=\sum_{i\in[d-1]}\bigg\langle{d\atop
i}\bigg\rangle\alpha^d(n-i).\nonumber
\end{align}
This computation of $\gamma^d(n)$ coincides with Kim's computation
\cite{kim} and we can also obtain this decomposition of
$\gamma^d(n)$ by a Worpitzky's result \cite{worpitzky}.

\subsubsection{The product of simplexes} For $j\in[l]$
let $\mathbf{0}^j$ be the zero vector and
$\mathbf{e}^j_0,\,\mathbf{e}^j_1,\ldots,\,\mathbf{e}^j_{d_j}$ be the
unit vectors in $\mathbbm{R}^{d_j+1}$. Writing
$$\alpha^{d_j}=conv\Big(\big\{\mathbf{e}^j_0,\mathbf{e}^j_1,\ldots,\mathbf{e}^j_{d_j}\big\}\Big)$$ for $j\in[l]$, we define
$$\alpha^{d_1,d_2,\ldots,d_l}=\underset{j\in[l]}{\prod}\alpha^{d_j}$$ to be the
\emph{standard product of simplexes}. Let
$d=\underset{j\in[l]}{\sum}d_j$ and
$S_{d_1,d_2,\ldots,d_l}$ be the set of permutations on the multiset
$\{1^{d_1},2^{d_2},\ldots,l^{d_l}\}$ where the number of $i$ in this
multiset is $d_i$. For a permutation $[a_1,a_2,\ldots,a_d]\in
S_{d_1,d_2,\ldots,d_l}$ denoting
$$\begin{cases}\big(A^1_0,A^2_0,\ldots,A^l_0\big)=(0,0,\ldots,0)\\
\big(A^1_i,A^2_i,\ldots,A^l_i\big)=\underset{k\in[i]}{\sum}
\mathbf{e}_{a_i}\end{cases}$$ for $i\in[d]$, we define
$\alpha^d[a_1,a_2,\ldots,a_d]$ to be the simplex with the vertex set
$$vert\big(
\alpha^d[a_1,a_2,\ldots,a_d]\big)=\bigg\{\prod_{j\in[l]}\mathbf{e}^{j}_{A^j_i}\,\bigg|\,i\in[d]_0\bigg\}$$ and
$$\mathcal{S}(\alpha^{d_1,d_2,\ldots,d_l})=\big\{\alpha^{d_1,d_2,\ldots,d_l}[a_1,a_2,\ldots,a_d]\,
\big|\,[a_1,a_2,\ldots,a_d]\in S_{d_1,d_2,\ldots,d_l}\big\}.$$

We claim that
\begin{equation}\label{dproduct}
\alpha^{d_1,d_2,\ldots,d_l}=\bigcup_{\alpha^d\in\mathcal{S}(\alpha^{d_1,d_2,\ldots,d_l})}\alpha^{d}.
\end{equation}
By definition
$$vert\Big(\alpha^{d_1,d_2,\ldots,d_l}\Big)\supseteq\underset{
\alpha^d\in\mathcal{S}(\alpha^{d_1,d_2,\ldots,d_l})}{\bigcup}vert(\alpha^d~),$$
thus
$$\alpha^{d_1,d_2,\ldots,d_l}
\supseteq\underset{\alpha^d\in
\mathcal{S}(\alpha^{d_1,d_2,\ldots,d_l})}{\bigcup}\alpha^d.$$
Consider the opposite inclusion $\subseteq$. Let
$$\mathbf{x}=\underset{j\in[l]}{\prod}
\big(x^j_0,x^j_1,\ldots,x^j_{d_j}\big)^*\in\alpha^{d_1,d_2,\ldots,d_l}.$$
Denoting $X^j_{k_j}=\underset{i\in[k_j]_0}{\sum}x^j_i$ for $(j,k_j)\in[l]\times [d_j]$, we let $\prec$ be the total order on
$X^j_{k_j}$ defined by $X^{j_1}_{k_{j_1}}\prec X^{j_2}_{k_{j_2}}$ if
one of the following is true:
$$\begin{cases}
X^{j_1}_{k_{j_1}}<X^{j_2}_{k_{j_2}}\\
X^{j_1}_{k_{j_1}}=X^{j_2}_{k_{j_2}}\ \text{for}\ j_1<j_2\\
X^{j_1}_{k_{j_1}}=X^{j_2}_{k_{j_2}}\ \text{for}\
j_1=j_2,\,k_{j_1}<k_{j_2}\end{cases}.$$ Assuming that
$X^{j_1}_{k_{j_1}}\prec X^{j_2}_{k_{j_2}}\prec\cdots \prec
X^{j_{d}}_{k_{j_{d}}}$, for $i\in[d]$ we define the $l$-tuples
$$\big(A^1_i,A^2_i,\ldots,A^l_i\big)=\sum_{k\in[i]}\mathbf{e}_{j_k}.$$ If we denote
$$\begin{cases}
X^0_{0}=\underset{j\in[l]}{\sum} x^j_0\\
X^{j_0}_{k_{j_0}}=0\end{cases},$$ then
$$\mathbf{x}=\sum_{i\in[d]}\Big(X^{j_i}_{k_{j_i}}-X^{j_{i-1}}_{k_{j_{i-1}}}\Big)
\prod_{j\in[l]}\mathbf{e}^j_{A^j_i}+X^0_0\prod_{j\in[l]}\mathbf{e}^j_0.
$$
The equation
$$\sum_{i\in[d+1]}\Big(X^{j_i}_{k_{j_i}}-X^{j_{i-1}}_{k_{j_{i-1}}}\Big)
+X^0_0=\sum_{j\in[l]}\sum_{i\in[d_j]_0}x^j_i=1$$
provides $\mathbf{x}\in\alpha^{d}[j_1,j_2,\ldots,j_d]$, thus
$$\prod_{j\in[l]}\alpha^{d_j}\subseteq
\underset{\alpha^d\in\mathcal{S}(\alpha^{d_1,d_2,\ldots,d_l})}{\bigcup}\alpha^d.$$

From the identity (\ref{dproduct}), we define a polytopal complex
$$\mathcal{C}_{d_1,d_2,\ldots,d_l}=\bigoplus_{\alpha^d\in
S(\alpha^{d_1,d_2,\ldots,d_l})}\alpha^d.$$ Then, by a simple
reasoning, $\mathcal{C}_{d_1,d_2,\ldots,d_l}$ is a pointed
triangulation of $\underset{j\in[l]}{\prod}\alpha^{d_j}$. We call
$\mathcal{C}_{d_1,d_2,\ldots,d_l}$ the \emph{standard pointed
triangulation} of $\underset{j\in[l]}{\prod}\alpha^{d_j}$.

For $i\in[d-1]_0$ let $\mathcal{S}_i(\alpha^{d_1,d_2,\ldots,d_l})$
be the set of $\alpha^d\big[a_1,a_2,\ldots,a_d\big]$ in
$\mathcal{S}_i(\alpha^{d_1,d_2,\ldots,d_l})$ such that for $k\in[d-1]_0$ the number of $k$ satisfying $a_k>a_{k+1}$ is $i$. We
define $$\bigg\langle{d_1,d_2,\ldots,d_l\atop
i}\bigg\rangle=\Big|\mathcal{S}_i(\alpha^{d_1,d_2,\ldots,d_l})\Big|.$$
The sets $\mathcal{S}_i(\alpha^{d_1,d_2,\ldots,d_l})$ satisfy
$$\mathcal{S}(\alpha^{d_1,d_2,\ldots,d_l})=\biguplus_{i\in[d-1]}
\mathcal{S}_i(\alpha^{d_1,d_2,\ldots,d_l}).$$ We say that an
$\alpha^d\big[a^1_1,a^1_2,\ldots,a^1_d\big]$ is \emph{adjacent} to
another $\alpha^d\big[a^2_1,a^2_2,\ldots,a^2_d\big]$ if
$$\big[a_1^2,a_2^2,\ldots,a_d^2\big]
=\big[a_1^1,\ldots,a_{k-1}^1,a_{k+1}^1,a_{k}^1,a_{k+2}^1,\ldots,a_d^1\big]$$
where $a^1_k<a^1_{k+1}$. By the definition of adjacency
$$dim\Big(\alpha^d\big[a_1^1,a_2^1,\ldots,\,a_{d}^1\big]\cap
\alpha^d\big[ a_1^2,a_2^2,\ldots,a_d^2\big]\Big)=d-1,$$ and if
$\alpha^d[
a_1^1,a_2^1,\ldots,a_{d}^1]\in\mathcal{S}_i(\alpha^{d_1,d_2,\ldots,d_l})$
then $$\alpha^d\big[
a_1^2,a_2^2,\ldots,a_{d}^2\big]\in\mathcal{S}_{i+1}(\alpha^{d_1,d_2,\ldots,d_l}).$$
In addition, if $\alpha^d[
a_1,a_2,\ldots,a_{d}]\in\mathcal{S}_i(\alpha^{d_1,d_2,\ldots,d_l})$
then the number of $d$-simplexes in
$\mathcal{S}(\alpha^{d_1,d_2,\ldots,d_l})$ that are adjacent to
$\alpha^d\big[ a_1,a_2,\ldots,a_{d}\big]$ is $i$.

For $\alpha^d\in\mathcal{S}(\alpha^{d_1,d_2,\ldots,d_l})$ let
$\mathcal{A}(\alpha^d)$ be the set of $d$-simplexes in
$\mathcal{S}(\alpha^{d_1,d_2,\ldots,d_l})$ that are adjacent to
$\alpha^d$ and let
$$C_{\alpha^d}(n)=S_{\alpha^d}(n)\setminus\bigcup_{\alpha_1^d\in
\mathcal{A}(\alpha^d)}S_{\alpha_1^d}(n).$$ Simply speaking,
$C_{\alpha^d}(n)$ is formed by successive facet-cuts on
$S_{\alpha^d}(n)$. By the geometric description of polytope numbers
$$S_{\alpha^{d_1,d_2,\ldots,d_l}}(n)=\biguplus_{\alpha_d\in
\mathcal{S}(\alpha^{d_1,d_2,\ldots,d_l})}C_{\alpha^d}(n),$$ and if
$\alpha^d\in\mathcal{S}_i(\alpha^{d_1,d_2,\ldots,d_l})$ then
successive $i$ facet-cuts on $S_{\alpha^d}$ yields
$$\big|C_{\alpha^d}(n)\big|=\alpha^d(n-i).$$
Therefore
\begin{align}
\alpha^{d_1,d_2,\ldots,d_l}(n)&=\big|S_{\alpha^{d_1,d_2,\ldots,d_l}}(n)\big|\label{dpsimplex}
\\
&=\bigg|\biguplus_{\alpha_d\in
\mathcal{S}(\alpha^{d_1,d_2,\ldots,d_l})}C_{\alpha^d}(n)\bigg|
=\bigg|\biguplus_{i\in[d-1]}\biguplus_{\alpha^d\in\mathcal{S}_i(\alpha^{d_1,d_2,\ldots,d_l})}C_{\alpha^d}(n)\bigg|
\nonumber\\
&=\sum_{i\in[d-1]}\sum_{\alpha^d\in\mathcal{S}_i(\alpha^{d_1,d_2,\ldots,d_l})}\big|C_{\alpha^d}(n)\big|
=\sum_{i\in[d-1]}\sum_{\alpha^d\in\mathcal{S}_i(\alpha^{d_1,d_2,\ldots,d_l})}\alpha^d(n-i)\nonumber\\
&=\sum_{i\in[d-1]}\bigg\langle{d_1,d_2,\ldots,d_l\atop
i}\bigg\rangle\alpha^d(n-i).\nonumber
\end{align}

\begin{remark} If $d_1=d_2=\cdots=d_l=1$, then
$\alpha^{d_1,d_2,\ldots,d_l}$ is an $l$-measure polytope. Thus the
identity (\ref{dpsimplex}) is a generalization of the identity 
(\ref{dmeasure}). In this point of view, we call the numbers $\Big\langle{d_1,d_2,\ldots,d_l\atop
i}\Big\rangle$ generalized Eulerian numbers.
\end{remark}

\subsection{Decomposition theorem 2 and Decomposition theorem 3}
Essentially, the same phenomenon describes Decomposition theorems 2
and 3. Therefore we apply these theorems at the same time.

\subsubsection{Simplex $\alpha^d$} For $\alpha_s^d$ let
$\mathbf{v}_{\alpha^d}=\mathbf{e}_0$. Suppose that
$\mathcal{C}_{\alpha_s^d}$ is the $V(\alpha_s^d)$-pointed triangulation.
Since a $k$-simplex in
$\mathcal{C}_{\alpha_s^d}$ that contains $\mathbf{v}_{\alpha^d}$ is
determined by $k$ vertexes among
$\mathbf{e}_1,\,\mathbf{e}_2,\ldots,\,\mathbf{e}_d$, the number of
such $k$-faces in $\mathcal{C}_{\alpha_s^d}$ is ${d\choose k}$.
Therefore, by Decomposition theorem 3-1,
$$\alpha^d(n)=\alpha_s^d(n)=\sum_{k\in[d]_0}{d\choose k}\alpha^k(n-k).$$

\subsubsection{Cross polytope $\beta^d$} For $\beta_s^d$ let
$\mathbf{v}_{\beta^d}=\mathbf{e}_d$. An $r$-simplex $\alpha^r$ in
$\mathcal{C}_{\beta_s^d}$ satisfies $\alpha^r\cap
int(\beta^d)\neq\emptyset$ if and only if
$\{\mathbf{e}_d,-\mathbf{e}_d\}\subseteq\alpha^r$, thus
$$\alpha^r=\{\mathbf{e}_d,-\mathbf{e}_d,a_{i_1}\mathbf{e}_{i_1},a_{i_2}\mathbf{e}_{i_2},\ldots,
a_{i_{r-1}}\mathbf{e}_{i_{r-1}}\}$$ where $a_{i_j}\in\{1,-1\}$. By
the way, the number of such $r$-simplexes in $\mathcal{C}_{\beta^d}$
is ${d-1\choose r-1}2^{r-1}$, thus Decomposition theorem 2 yields
$$\beta^d(n)=\beta_s^d(n)=\sum_{r\in[d]}(-1)^{d-r}{d-1\choose
r-1}2^{r-1}\alpha^r(n).$$ This computation of $\beta^d(n)$ coincides
with Kim's computation of $\beta^d(n)$ \cite{kim}. Similarly,
Decomposition theorem 3-2 provides
$$\beta^d(n)=\sum_{k\in[d]}{d-1\choose
k-1}2^{k-1}\alpha^k(n-k+1).$$ In addition, if we consider the number
of $k$-simplexes in $\mathcal{C}_{\beta_s^d}$ that are on the
boundary of $\beta_s^d$ and contain $\mathbf{v}_{\beta_s^d}$, which
is $2^k{d-1\choose k}$ by a simple computation, then, by
Decomposition theorem 3-1,
$$\beta^d(n)=\sum_{k\in[d]_0}\bigg\{2^k{d-1\choose
k}+2^{k-1}{d-1\choose k-1}\bigg\}\alpha^k(n-k).$$

\subsubsection{Measure polytope $\gamma^d$} For $\gamma_s^d$ let
$\mathbf{v}_{\gamma^d}=\mathbf{0}$. If $\alpha^r$ is an $r$-simplex
in $\mathcal{C}_{\gamma^d}$ satisfying $\alpha^r\cap
int(\gamma^d)\neq\emptyset$, then
$\{\mathbf{0},\mathbf{1}\}\subseteq\alpha^r$. Moreover, a
$d$-simplex $\alpha^d[a_1,a_2,\ldots,a_d]\in\mathcal{C}_{\gamma^d}$
satisfying $\alpha^r\subseteq\alpha^d[a_1,a_2,\ldots,a_d]$ exists,
thus the vertexes of $\alpha^r$ are
$$\mathbf{0},\,\sum_{k_1\in[d_1]}\mathbf{e}_{a_{k_1}},\,\sum_{k_2\in[d_2]}\mathbf{e}_{a_{k_2}},\ldots,\,
\sum_{k_{r-1}\in[d_{r-1}]}\mathbf{e}_{a_{k_{r-1}}},\,\mathbf{1}$$
where $1\leq d_1<d_2<\cdots<d_{r-1}\leq d-1$. Consequently, the
vertexes of $\alpha^r$ correspond to an ordered partition of
$[d]$ into $r$ sets. The number of such partitions is
$r!S(d,r)$ where $S(d,r)$ is the Stirling number of the second kind,
which is the number of ways to partition a $d$-set into $r$ sets. If
we assume that
$$a_1>a_2>\cdots>a_{d_1},a_{d_1+1}>a_{d_1+2}>\cdots>a_{d_2},\ldots,
a_{d_{r-1}+1}>a_{d_{r-1}+2}>\cdots>a_d,$$ that is,
$[a_1,a_2,\ldots,a_d]$ has at least $d-r$ descents, then the number
of such $\alpha^r$ is $$\sum_{i\in[d-1]\setminus[d-1-r]}\bigg\langle{d \atop
i}\bigg\rangle{i\choose d-r}.$$ Therefore Decomposition theorem 2
yields
\begin{align*}
\gamma^d(n)&=\gamma_s^d(n)=\sum_{r\in[d]}(-1)^{d-r}r!S(d,r)\alpha^r(n)\\
&=\sum_{r\in[d]}(-1)^{d-r}\bigg\{\sum_{i\in [d-1]\setminus [d-1-r]}\bigg\langle{d
\atop i}\bigg\rangle{i\choose d-r}\bigg\}\alpha^r(n). \end{align*}
Similarly, by Decomposition theorem 3-2,
\begin{align*}
\gamma^d(n)&=\sum_{r\in[d]}r!S(d,r)\alpha^r(n-r+1)\\
&=\sum_{r\in[d]}\bigg\{\sum_{i\in[d-1]\setminus[d-1-r]}\bigg\langle{d \atop
i}\bigg\rangle{i\choose d-r}\bigg\}\alpha^r(n-r+1).
\end{align*}

Let $\alpha^r$ be an $r$-simplex in $\mathcal{C}_{\gamma_s^d}$ that
contains $\mathbf{0}$ and is a subpolytope of $\alpha^d[a_1,a_2,\ldots,a_d]$. Then
$$\alpha^r=vert\bigg(\bigg\{
\mathbf{0},\,\sum_{k_1\in[d_1]}\mathbf{e}_{a_{k_1}},\sum_{k_2\in[d_2]}\mathbf{e}_{a_{k_2}},\ldots,
\sum_{k_{r-1}\in[d_{r-1}]}\mathbf{e}_{a_{k_{r-1}}}\bigg\}\bigg)$$
where $1\leq d_1<d_2<\cdots<d_{r-1}\leq d$. Thus to choose an
$r$-simplex in $\mathcal{C}_{\gamma_s^d}$ is equivalent to construct
an ordered partition of a subset of $[d]$ with $j$
elements into $r$ sets where $j\in[d]\setminus[r-1]$. It follows that the
number of such $r$-simplexes in $\mathcal{C}_{\gamma_s^d}$ is
$$\sum_{j\in[d]\setminus[r-1]}{d\choose j}r!S(j,r).$$
Decomposition theorem 3-1 yields
$$\gamma^d(n)=\sum_{r\in[d]_0}\bigg(\sum_{j\in[d]\setminus[r-1]}{d\choose
j}r!S(j,r)\bigg)\alpha^r(n-r).$$

\subsubsection{The product of simplexes} For
$\alpha^{d_1,d_2,\ldots,d_l}$ let
$$\mathbf{v}_{\alpha^{d_1,d_2,\ldots,d_l}}=\prod_{j\in[l]}\mathbf{e}^j_0.$$
Suppose that $\alpha^r\subseteq\alpha^d[a_1,a_2,\ldots,a_d]$ and for $i_1<i_2<\cdots<i_r$ let
$$\prod_{j\in[l]}\mathbf{e}^j_0,\prod_{j\in[l]}\mathbf{e}^j_{A^j_{i_1}},
\prod_{j\in[l]}\mathbf{e}^j_{A^j_{i_2}},\ldots,\,\prod_{j\in[l]}\mathbf{e}^j_{A^j_{i_r}
}$$ be the vertexes of $\alpha^r$. The
simplex $\alpha^r$ satisfies $\alpha^r\cap
int(\alpha^{d_1,d_2,\ldots,d_l})\neq\emptyset$ if and only if there
is a point
$\mathbf{x}=\underset{j\in[l]}{\prod}\big(x^j_0,x^j_1,\ldots,x^j_{d_j}\big)^*$
in $\alpha^r$ such that every entry of $\mathbf{x}$ is nonzero,
equivalently,
\begin{equation}\label{condition1}
\Big\{0,A^j_{i_1},A^j_{i_2},\ldots,A^j_{i_r}\Big\}=[d_j]_0
\end{equation}
for $j\in[l]$. Therefore, if we denote
$\big(A^1_{i_0},A^2_{i_0},\ldots,A^l_{i_0}\big)=\mathbf{0}$, then
the number of $r$-simplexes $\alpha^r$ satisfying $\alpha^r\cap
int(\alpha^{d_1,d_2,\ldots,d_l})\neq\emptyset$ is the number of elements in the set $\mathcal{A}$ composed of
$$\Big(\big(A^1_{i_1},A^2_{i_1},\ldots,A^l_{i_1}\big),\big(A^1_{i_2},A^2_{i_2},\ldots,A^l_{i_2}\big),
\ldots,\big(A^1_{i_r},A^2_{i_r},\ldots,A^l_{i_r}\big)\Big)$$
that satisfies
\begin{equation}\label{condition2}
\left\{\begin{array}{ll}
\big(A^1_{i_m},A^2_{i_m},\ldots,A^l_{i_m}\big)-
\big(A^1_{i_{m-1}},A^2_{i_{m-1}},\ldots,A^l_{i_{m-1}}\big)>
\mathbf{0}&\text{for}\ m\in[r]\\
\Big\{0,A^j_{i_1},A^j_{i_2},\ldots,A^j_{i_r}\Big\}=[d_j]_0&
\text{for}\ j\in[l]\end{array}\right..\end{equation}

For $k\in[r]$ let $\mathcal{A}_k$ be the subset of $\mathcal{A}$ each of whose
elements satisfies
\begin{equation*}
\big(A^1_{i_k},A^2_{i_k},\ldots,A^l_{i_k}\big)-
\big(A^1_{i_{k-1}},A^2_{i_{k-1}},\ldots,A^l_{i_{k-1}}\big)
=\mathbf{0}.
\end{equation*}
The number of
$$\Big(\big(A^1_{i_1},A^2_{i_1},\ldots,A^l_{i_1}\big),\big(A^1_{i_2},A^2_{i_2},\ldots,A^l_{i_2}\big),
\ldots,\big(A^1_{i_r},A^2_{i_r},\ldots,A^l_{i_r}\big)\Big)$$
satisfying the condition (\ref{condition2}) is the number of
elements in
$\mathcal{A}\setminus\Big(\underset{k\in[r]}{\bigcup}\mathcal{A}_k\Big)$. By the
principle of inclusion and exclusion \cite{stanley}
$$\bigg|\mathcal{A}\setminus\bigg(\bigcup_{k\in[r]}\mathcal{A}_k\bigg)\bigg|
=\sum_{k\in[r]_0}(-1)^k{r\choose k}\prod_{j\in[l]}{r-k\choose d_j},$$ thus
Decomposition theorems 2 and 3-2 yield
$$\begin{cases}
\alpha^{d_1,d_2,\ldots,d_l}(n)=
\underset{r\in[d]}{\sum}(-1)^{d-r}\bigg\{\underset{k\in[r]_0}{\sum}(-1)^k{r\choose
k}\underset{j\in[l]}{\prod}{r-k\choose d_j}\bigg\}\alpha^r(n)\\
\alpha^{d_1,d_2,\ldots,d_l}(n)=\underset{r\in[d]}{\sum}\bigg\{\underset{k\in[r]_0}{\sum}(-1)^k{r\choose
k}\underset{j\in[l]}{\prod}{r-k\choose d_j}\bigg\}\alpha^r(n-r+1)\end{cases},$$
respectively.

 We can similarly compute the number of $r$-simplexes
in $\mathcal{C}_{d_1,d_2,\ldots,d_l}$ that contain
$\mathbf{v}_{\alpha^{d_1,d_2,\ldots,d_l}}$. It is the number of
$$\Big(\big(A^1_{i_1},A^2_{i_1},\ldots,A^l_{i_1}\big),\big(A^1_{i_2},A^2_{i_2},\ldots,A^l_{i_2}\big),
\ldots,\big(A^1_{i_r},A^2_{i_r},\ldots,A^l_{i_r}\big)\Big)$$ obeying
\begin{equation*}
\left\{\begin{array}{ll}
\big(A^1_{i_m},A^2_{i_m},\ldots,A^l_{i_m}\big)-
\big(A^1_{i_{m-1}},A^2_{i_{m-1}},\ldots,A^l_{i_{m-1}}\big)>
\mathbf{0}& \text{for}\ m\in[r]\\
\Big\{0,A^j_{i_1},A^j_{i_2},\ldots,A^j_{i_r}\Big\}=[d^1_j]_0&
\text{for}\ j\in[l]\ \text{and}\ d^1_j\in[d_j]_0\end{array}\right..\end{equation*}
By the principle of inclusion
and exclusion, this number is
$$\sum_{k\in[r]_0}(-1)^k{r\choose k}\prod_{j\in[l]}{d_j+r-k\choose
   r-k},$$
therefore Decomposition theorem 3-1 furnishes
$$\alpha^{d_1,d_2,\ldots,d_l}(n)=\sum_{r\in[d]_0}\bigg\{\sum_{k\in[r]_0}(-1)^k{r\choose
k}\prod_{j\in[l]}{d_j+r-k\choose r-k}\bigg\}\alpha^r(n-r).$$

\section{Applications of polytopes numbers\label{section6}}

Decomposition theorems are methods to represent polytope numbers by
sums of simplex numbers and the vertex description of polytope
numbers gives a relation between polytope numbers and a set of
chains in posets formed by faces of polytopes. Using these facts, we
consider applications of polytope numbers to several research topics.
These topics are composed of generalized Eulerian numbers, lattice
paths, plane partitions, and Young tableaux. For two polytopes $P_1$
and $P_2$ we say that $P_1$ is a vertex subpolytope of $P_2$ if
$vert(P_1)\subseteq vert(P_2)$.

\subsection{Generalized Eulerian numbers\label{subsection6.1}}
In this subsection, we denote the unit vectors of $\mathbb{R}^{d+1}$
by $\mathbf{e}_0,\mathbf{e}_1,\ldots,\mathbf{e}_d$ and we define
$$\mathbf{e}_{i_1,i_2,\ldots,i_l}=\mathbf{e}_{i_1}\times\mathbf{e}_{i_2}\times
\cdots\times\mathbf{e}_{i_l}$$ where $(i_1,i_2,\ldots,i_l)\in[d_1]_0\times[d_2]_0\times\cdots[d_l]_0$.

Letting $d=\underset{i\in[l]}{\sum}d_i$, we define
$$L(i_1,i_2,\ldots,i_{d+1})=\{\mathbf{e}_{i_1},\mathbf{e}_{i_1}+\mathbf{e}_{i_2},\ldots,
\mathbf{e}_{i_1}+\mathbf{e}_{i_2}+\cdots+\mathbf{e}_{i_{d+1}}\}$$
to be the lattice path
$$\mathbf{0}\rightarrow\mathbf{e}_{i_1}\rightarrow\mathbf{e}_{i_1}+\mathbf{e}_{i_2}\rightarrow
\cdots\rightarrow\mathbf{e}_{i_1}+\mathbf{e}_{i_2}+\cdots+\mathbf{e}_{i_d}$$
from $(0,0,\ldots,0)$ to $(d_1,d_2,\ldots,d_l)$,
and $\mathcal{L}(d_1,d_2,\ldots,d_l)$ to be the set of such $L(i_1,i_2,\ldots,i_{d+1})$.
For $L(i_1,i_2,\ldots,i_{d+1})\in\mathcal{L}(d_1,d_2,\ldots,d_l)$ let
$$\alpha^d(i_1,i_2,\ldots,i_{d+1})=conv\Big(
\big\{\mathbf{e}_{j_1,j_2,\ldots,j_l}\,\big|\,(j_1,j_2,\ldots,j_l)\in
L(i_1,i_2,\ldots,i_{d+1})\big\}\Big).$$ be a $d$-simplex, which is a
vertex subpolytope of $\alpha^{d_1,d_2,\ldots,d_l}$.

For a lattice path $L(i_1,i_2,\ldots,i_{d+1})$ in
$\mathcal{L}(d_1,d_2,\ldots,d_l)$, we define a descent of
$L(i_1,i_2,\ldots,i_{d+1})$ to be an index $p$ such that
$i_p>i_{p+1}$. Then we can easily show that
\begin{equation}\label{edecomposition}
\prod_{i\in[l]}\alpha^{d_i}=\underset{L(i_1,i_2,\ldots,i_{d+1})\in
\mathcal{L}(d_1,,d_2,\ldots,d_l)}{\bigcup}\alpha^d(i_1,i_2,\ldots,i_{d+1}).\end{equation}
By Decomposition theorem 1,
$$\prod_{i\in[l]}\alpha^{d_i}(n)=\sum_{i\in[d-1]_0}a_i\alpha^d(n-i)$$ where $a_i$ is the number
of lattice paths in $\mathcal{L}(d_1,d_2,\ldots,d_l)$ with $i$
descents. Since $a_i=\Big\langle {l\atop i}\Big\rangle$ for
$d_1=d_2=\cdots=d_l=2$, the numbers $a_0,a_1,\ldots,a_{d-1}$ are
\emph{generalized Eulerian numbers}. By
$$\begin{pmatrix}
\underset{i\in[l]}{\prod}\alpha^{d_i}(1)\\
\underset{i\in[l]}{\prod}\alpha^{d_i}(2)\\
\cdots\\
\underset{i\in[l]}{\prod}\alpha^{d_i}(d)
\end{pmatrix}
=\begin{pmatrix} \alpha^d(1)&0&\cdots&0\\
\alpha^d(2)&\alpha^d(1)&\cdots&0\\
\cdots&\cdots&\cdots&\cdots\\
\alpha^d(d)&\alpha^d(d-1)&\cdots&\alpha^d(1)
\end{pmatrix}\begin{pmatrix}
a_{0}\\a_{1}\\ \cdots\\a_{d-1}
\end{pmatrix},
$$
the numbers $a_0,a_1,\ldots,a_{d-1}$ are
$$\begin{pmatrix}
a_{0}\\a_{1}\\ \cdots\\a_{d-1}
\end{pmatrix}
=
\begin{pmatrix}
(-1)^0{d+1\choose 0}&0&\cdots&0\\
(-1)^1{d+1\choose 1}&(-1)^0{d+1\choose 0}&\cdots&0\\
\cdots&\cdots&\cdots&\cdots\\
(-1)^{d-1}{d+1\choose d-1}&(-1)^{d-2}{d+1\choose
d-2}&\cdots&(-1)^0{d+1\choose 0}
\end{pmatrix}
\begin{pmatrix}
\underset{i\in[l]}{\prod}\alpha^{d_i}(1)\\
\underset{i\in[l]}{\prod}\alpha^{d_i}(2)\\
\cdots\\
\underset{i\in[l]}{\prod}\alpha^{d_i}(d)
\end{pmatrix}.
$$
Therefore
$$a_i=\sum_{j\in[i]}(-1)^{j}{d+1\choose
j}\prod_{k\in[l]}\alpha^{d_k}(i+1-j).$$

\subsection{Lattice paths and plane partitions}
Let $P_L$ be a $d$-dimensional vertex subpolytope of
$\alpha^{d_1,d_2,\ldots,d_l}$ with $d_1\ge d_2\ge\cdots\ge d_l$ each of whose verteices  $\mathbf{e}_{i_1,i_2,\ldots,i_l}$ satisfies
$i_1\ge\i_2\ge\cdots\ge i_l$. If we use the decomposition
(\ref{edecomposition}), then the $d$-simplexes
$\alpha^d(i_1,i_2,\ldots,i_{d+l})$ in the decomposition
(\ref{edecomposition}) correspond to the lattice paths from
$(0,0,\ldots,0)$ to $(d_1,d_2,\ldots,d_l)$ such that every point
$(x_1,x_2,\ldots,x_d)$ in these lattice paths satisfies $x_1\ge
x_2\ge\cdots\ge x_l$. By Decomposition theorem 1,
$$P_L(n)=\sum_{i\in[d-1]}a_{i}\alpha^d(n-i)$$
where $a_i$ is the number of simplexes
$\alpha^d(i_1,i_2,\ldots,i_{d+1})$ such that the number of descents in
$L(i_1,i_2,\ldots,i_{d+1})$ is $i$ and every point
$(x_1,x_2,\ldots,x_l)$ in the lattice path $L(i_1,i_2,\ldots,i_{d+1})$
satisfies $x_1\ge x_2\ge\cdots\ge x_l$. Krattenthaler computed the
number $\underset{i\in[d-1]_0}{\sum}a_i$ \cite{krattenthaler} and $a_i$ are
refinements of $\underset{i\in[d-1]_0}{\sum}a_i$.

We now compute the numbers $a_i$. According to the decomposition
(\ref{edecomposition}), we can define a partial order of the
vertexes in $P_L$ by
$$\mathbf{e}_{a_{11},a_{12},\ldots,a_{1l}}\ge
\mathbf{e}_{a_{21},a_{22},\ldots,a_{2l}}$$ if $$a_{11}\leq
a_{21},a_{12}\leq a_{22},\ldots,a_{1l}\leq a_{2l}.$$ Therefore we
need to consider the number of different sums of $n-1$ vertexes
$\mathbf{e}_{a_{i1},a_{i2},\ldots,a_{il}}$ for $i\in[n-1]$ in
$P_L$ such that
$$\mathbf{e}_{a_{i1},a_{i2},\ldots,a_{il}}\ge\mathbf{e}_{a_{i+1,1},a_{i+1,2},\ldots,a_{i+1,l}}$$
to use the vertex description of polytope numbers. By a simple
computation, $P_L(n)$ is the number of plane partitions with entries
$a_{ij}$ satisfying $a_{ij}\in[d_1]$. Since the number of such
plane partitions is
$$P(n-1,l,d_1)=\prod_{i\in[n-1]}\prod_{j\in[l]}\prod_{k\in[d_1]}\frac{i+j+k-1}{i+j+k-2},$$
where $P(0,l,d_1)=1$ \cite{macmahon}, the vertex
description of polytope numbers for $P_L$ yields
$$P_L(n)=P(n-1,1,d_1).$$
The method to compute coefficients in the decomposition form of Decomposition
theorem 1 yields
$$\begin{pmatrix} a_{0}\\a_{1}\\
\cdots\\a_{d-1}
\end{pmatrix}
=
\begin{pmatrix}
(-1)^0{d+1\choose 0}&0&\cdots&0\\
(-1)^1{d+1\choose 1}&(-1)^0{d+1\choose 0}&\cdots&0\\
\cdots&\cdots&\cdots&\cdots\\
(-1)^{d-1}{d+1\choose d-1}&(-1)^{d-2}{d+1\choose
d-2}&\cdots&(-1)^0{d+1\choose 0}
\end{pmatrix}
\begin{pmatrix}
P(0,l,d_1)\\P(1,l,d_1)\\ \cdots\\P(d-1,l,d_1)
\end{pmatrix}.$$
Therefore
$$a_i=\sum_{j\in[i]_0}(-1)^j{d+1\choose j}P(i-j,l,d_1).$$

If we let $l=2$, then $a_i$ is the number of lattice paths from
$(0,0)$ to $(d_1,d_2)$ each of whose unit paths is either east or
north and the first unit path is east, which never crosses the line
$x_1=x_2$, and whose number of consecutive north-east paths is $i$.
The numbers $a_i$ are refinements of the Lobb number
$\underset{i\in[d-1]_0}{\sum}a_i$ \cite{lobb}. In particular, if $d_1=d_2$
then $a_{i}=N(d_1,i+1)$ is the Narayana number \cite{macmahon,
narayana}.

If we let $d_1=d_2=\cdots=d_l$, then $a_i$ are \emph{higher
dimensional Narayana numbers} \cite{sulanke}.

\subsection{Young tableaux and plane partitions}
Let $P_Y$ be a vertex subpolytope of
$\gamma^{l\cdot m}=\underset{i\in[l\cdot m]}{\prod}conv\big(\{\mathbf{e}_0,\mathbf{e}_1\}\big)$
whose vertexes are
$$\prod_{(i,j,a_{ij})\in S}\mathbf{e}_{a_{ij}}$$ where
$S$ is a subset of $[l]\times[m]\times[1]_0$ such that
each pair of $(i_1,j_1,a_{i_1j_1})$ and $(i_2,j_2,a_{i_2j_2})$ in $S$
satisfies that $$a_{i_1j_1}\ge a_{i_2j_2}\ \text{if and only if}\
i_1\leq i_2\ \text{and}\ j_1\leq j_2.$$ Then $P_Y$ is an $(l\cdot
m)$-polytope and there is a one-to-one correspondence between the
vertexes of $P_Y$ and the partitions whose size of largest part is
at most $m$ and the number of parts is at most $l$. Note that a
partition $\lambda$ is a finite sequence of positive integers
$\lambda_1,\lambda_2,\ldots,\lambda_n$ satisfying
$$\lambda_1\ge\lambda_2\ge\cdots\ge\lambda_n$$ and we call $\lambda_1,\lambda_2,\ldots,\lambda_n$
the parts of $\lambda$ and $n$ the number of parts in $\lambda$.
Using these two properties of the polytope $P_Y$, we can represent
every vertex of $P_Y$ as follows: Let $(s_1,s_2,\ldots,s_l)$ be a
partition whose parts are $s_1,s_2,\ldots,\,s_l$ with $m\ge s_1\ge
s_2\ge\ldots\ge s_l$. We define
$$\mathbf{v}(s_1,s_2,\ldots,s_l)=\prod_{(i,j,a_{ij})\in S}\mathbf{e}_{a_{ij}}$$
where
$$a_{ij}=\begin{cases}
1& \text{for}\ j\in[s_i]\\
0& \text{for}\ j\in[m]\setminus[s_i]
\end{cases}$$ when $i\in[l]$.
We assign a partial order to $vert(P_Y)$ by
$$\mathbf{v}\big(s^1_1,s^1_2,\ldots,s^1_l\big)\ge
\mathbf{v}(s^2_1,s^2_2,\ldots,s^2_l)\ \text{if}\ s^1_1\leq
s^2_1,s^1_2\leq s^2_2,\ldots,s^1_l\leq s^2_l.$$

Let $\mathcal{Y}(P_Y)$ be the set of
$$\Big\{\mathbf{v}(s^0_1,s^0_2,\ldots,s^0_l),\mathbf{v}(s^1_1,s^1_2,\ldots,s^1_l),\ldots,
\mathbf{v}\big(s^{l\cdot m}_1,s^{l\cdot m}_2,\ldots,s^{l\cdot
m}_l\big)\Big\}\subseteq vert(P_Y)$$ satisfying
$$\begin{cases}
\big(s^0_1,s^0_2,\ldots,s^0_l\big)=(0,0,\ldots,0)\\
\big(s^{i}_1,s^{i}_2,\ldots,s^{i}_l\big)-
\big(s^{i-1}_1,s^{i-1}_2,\ldots,s^{i-1}_l\big)=\mathbf{e}_j\ \text{for}\
i\in[l\cdot m]\end{cases}$$ for some $j\in[l]$.
For $Y\in\mathcal{Y}(P_Y)$ we define an $(l\cdot
m)$-simplex
$$\alpha^{l\cdot m}_Y=conv(Y).$$ Since $P_Y$ is a
vertex subpolytope of the product of simplexes, decompositions of the
product of simplexes gives
$$P_Y=\bigcup_{Y\in\mathcal{Y}(P_Y)}\alpha^{l\cdot m}_Y.$$ Moreover,
$\underset{Y\in\mathcal{Y}(P_Y)}{\bigoplus}\alpha^{l\cdot m}_Y$ is a pointed
triangulation of $P_Y$, therefore the combination of these results
and the vertex description of polytope numbers yields
$$P_Y(n)=P(n-1,l,m).$$

For each $Y\in\mathcal{Y}(P_Y)$ we construct an $l\times m$ Young
tableau $T(Y)$ with entries in $[l\cdot m]$ whose entries are
strictly decreasing in each row and column as follows: For an
element $\mathbf{v}(s_1,s_2,\ldots,s_l)$ of $Y$ we define an
$l\times m$ matrix $M\big(\mathbf{v}(s_1,s_2,\ldots,s_l)\big)$ with
entries
$$M\big(\mathbf{v}(s_1,s_2,\ldots,s_l)\big)_{ij}=\begin{cases}
1& \text{for}\ j\in[s_i]\\
0& \text{for}\ j\in[m]\setminus[s_i]
\end{cases}$$ where $i\in[l]$.
The the $(i,j)$ entry of $T(Y)$ is the number of $\mathbf{v}\in Y$
such that $M(\mathbf{v})_{ij}=1$. Then $T(Y)$ is an $l\times m$
Young tableau with entries in $[l\cdot m]$ where each row and
column of $T(Y)$ has strictly decreasing entries.

Let
$$Y=\Big\{\mathbf{v}\big(s^0_1,s^0_2,\ldots,s^0_l\big),
\mathbf{v}\big(s^1_1,s^1_2,\ldots,s^1_l\big),\ldots,\mathbf{v}\big(s^{l\cdot
m}_1,s^{l\cdot m}_2,\ldots, s^{l\cdot m}_l\big)\Big\}$$ be an
element of $\mathcal{Y}(P)$ with
$$\begin{cases}
\big(s^0_1,s^0_2,\ldots,s^0_l\big)=(0,0,\ldots,0)\\
$$\big(s^i_1,s^i_2,\ldots,s^i_l\big)-\big(s^{i-1}_1,s^{i-1}_2,\ldots,s^{i-1}_l\big)=\mathbf{e}_j\
\text{for}\ i\in[l]\end{cases}$$ for some $j\in[l]$.
We define a \emph{descent} of $T(Y)$ to be an index $i\in[l\cdot m-1]$
such that if
$$\begin{cases}
\big(s^{i-1}_1,s^{i-1}_2,\ldots,s^{i-1}_l\big)-\big(s^{i}_1,s^{i}_2,\ldots,s^{i}_l\big)=\mathbf{e}_{j_1}\\
\big(s^{i}_1,s^{i}_2,\ldots,s^{i}_l\big)-\big(s^{i+1}_1,s^{i+1}_2,\ldots,s^{i+1}_l\big)=\mathbf{e}_{j_2}\end{cases}$$
then $j_1>j_2$. If we use both the result in decompositions of the
product of simplexes and that in Decomposition theorem 1, then the
coefficient $a_i$ in
$$P_Y(n)=\sum_{i\in[l\cdot m-1]_0}a_i\alpha^{l\cdot m}(n-i)$$
is the number of $l\times m$ Young tableaux that have exactly $i$
descents. The method of computing coefficients in
Decomposition theorem $1$ yields
$$\begin{pmatrix} a_{0}\\a_{1}\\
\cdots\\a_{l\cdot m-1}
\end{pmatrix}
=A
\begin{pmatrix}
P(0,l,m)\\P(1,l,m)\\ \cdots\\P(l\cdot m-1,l,m)
\end{pmatrix},$$
where
$$
A=\begin{pmatrix}
(-1)^0{l\cdot m+1\choose 0}&0&\cdots&0\\
(-1)^1{l\cdot m+1\choose 1}&(-1)^0{l\cdot m+1\choose 0}&\cdots&0\\
\cdots&\cdots&\cdots&\cdots\\
(-1)^{l\cdot m-1}{l\cdot m+1\choose l\cdot m-1}&(-1)^{l\cdot
m-2}{l\cdot m+1\choose l\cdot m-2}&\cdots&(-1)^0{l\cdot m+1\choose
0}
\end{pmatrix}.
$$
 Therefore
$$a_i=\sum_{j\in[i]_0}(-1)^j{l\cdot m+1\choose j}P(i-j,l,m).$$
The number of $l\times m$ Young tableaux with entries in
$[l\cdot m]$ is
$$(l\cdot m)!\frac{\underset{j\in[m-1]}{\prod}j!}{\underset{k\in[l+m-1]}{\prod}k!}$$
by the hook length formula \cite{frame}, thus
\begin{align*}\sum_{i\in[l\cdot m-1]_0}a_i&=\sum_{i\in[l\cdot
m-1]_0}\sum_{j\in[i]_0}(-1)^j{l\cdot m+1\choose j}P(i-j,l,m)\\
&=(l\cdot m)!\frac{\underset{j\in[m-1]}{\prod}j!}{\underset{k\in[l+m-1]}{\prod}k!}
.\end{align*}
\nocite{euler,gabrielov,gelfand,vanlint,lobb,neumann,schoute,stanley0,young}

\bibliographystyle{elsarticle-num}
\bibliography{Polytope_numbers_and_their_properties}

\end{document}